\documentclass[12pt]{amsart}
\usepackage{latexsym}
\usepackage{amsmath,amsfonts,amssymb,amsthm}
\usepackage{color}
\usepackage{xypic}
\input xy
\xyoption{all}
\newcommand{\kitem}{\begin{itemize}\vspace{-2ex}}
\newcommand{\kenditem}{\vspace{-1ex}\end{itemize}}
\newcommand{\A}{\mathbb A}
\newcommand{\C}{\mathbb C}

\newcommand{\PP}{\mathbb P}

\newcommand{\R}{\mathbb R}
\newcommand{\T}{\mathbb T}
\newcommand{\V}{\mathbb V}
\newcommand{\Z}{\mathbb Z}

\newcommand{\CN}{{\mathcal N}}

\DeclareMathOperator{\Spec}{Spec}

\DeclareMathOperator{\cone}{cone}

\DeclareMathOperator{\im}{im}

\DeclareMathOperator{\Cone}{Cone}

\DeclareMathOperator{\spann}{span}

\newcommand{\conv}{\operatorname{conv}} %
\DeclareMathOperator{\innt}{int}

\DeclareMathOperator{\Pic}{Pic}

\newcommand{\gHom}{\mbox{\rm Hom}}

\reversemarginpar
\newcommand{\kkk}[1]{}
\newcommand{\ktrash}[1]{{}}

\newcommand{\veee}{{\scriptscriptstyle\vee}}
\newcommand{\kb}{{\kss \bullet}}
\newcommand{\ku}{\underline}
\newcommand{\ko}{\overline}
\newcommand{\ks}{\scriptstyle}

\newcommand{\kf}{\footnotesize}
\newcommand{\kss}{\scriptscriptstyle}
\newcommand{\surj}{\rightarrow\hspace{-0.8em}\rightarrow}

\newcommand{\kP}{{\Delta}}

\newcommand{\kT}{{D}}

\newcommand{\kn}{{k}}

\newcommand{\kPV}[1]{{V^{#1}}}

\newcommand{\an}{{m}}
\newcommand{\kPF}{{F}}

\newcommand{\kPf}{{f}}

\newcommand{\kQ}{Q}
\newcommand{\kQP}{P}
\newcommand{\kQT}{T}
\newcommand{\wt}{\theta}
\newcommand{\wtc}{\theta^c}

\newcommand{\kA}{{\mbox{Inc}}}
\newcommand{\kH}{I\!\!H}
\newcommand{\kF}{I\!\!F}

\newcommand{\kD}{{\nabla}}
\newcommand{\kDQ}{{\kD(\kQ)}}

\newcommand{\kDQwt}{{\kD(\kQ,\wt)}}

\newcommand{\kDQwtc}{{\kD(\kQ,\wtc)}}
\newcommand{\kDPwt}{{\kD(\kQP,\wt)}}
\newcommand{\kDPwtc}{{\kD(\kQP,\wtc_\kQ)}}
\newcommand{\kDst}{{\Delta}}
\newcommand{\kDstQ}{{\kDst(\kQ)}}

\newcommand{\kwtTrees}{{{\mathcal T}^\wt}}

\newcommand{\tail}{t}
\newcommand{\head}{h}

\newcommand{\kG}{\Gamma}

\newcommand{\kGdbl}{{\kG^{\mbox{\tiny dbl}}}}
\newcommand{\kGopp}{{\kG^{\mbox{\tiny opp}}}}

\newcommand{\fence}{fence}

\newtheorem{theorem}{Theorem}

\newtheorem{lemma}[theorem]{Lemma}

\newtheorem{corollary}[theorem]{Corollary}

\newtheorem{proposition}[theorem]{Proposition}

\theoremstyle{definition}
\newtheorem{definition}[theorem]{Definition}

\newtheorem{example}[theorem]{Example}

\theoremstyle{remark}
\newtheorem*{remark}{Remark}

\numberwithin{equation}{section}

\newcommand{\toric}[1]{\T\V({#1})}  %
\newcommand{\normal}{\CN}  %
\newcommand{\klabel}[1]{\label{#1}\kkk{#1}}
\definecolor{changes}{rgb}{0.1,0.65,0.03}

\newcommand{\kst}{\,|\;}
\newcommand{\kSt}{\,\big|\;}

\newcounter{Abschnitt}[section]
\newcommand{\neu}[1]{\protect\refstepcounter{Abschnitt}\protect
   \label{#1}\vspace{1ex}
   {\bf (\protect\arabic{section}.\protect\arabic{Abschnitt})}
                     $\qquad$\kkk{#1}}
\newcommand{\zitat}[2]{(\protect\ref{#1}.\protect\ref{#1-#2})}

\begin{document}
\title{Smoothing of Quiver Varieties}
\author{Klaus Altmann \quad\qquad Duco van Straten }
\date{}
\maketitle
\begin{abstract}
We show that Gorenstein singularities
that are cones over singular Fano varieties
provided by so-called flag quivers are smoothable in codimension three.
Moreover, we give a precise characterization about the
smoothability in codimension three of the Fano variety itself.
\end{abstract}

\section{Introduction}\label{intro}

\neu{intro-setup}
Quivers and varieties of quiver representations appear in various places in
mathematics, see for example \cite{tilt}. In \cite{flags} it was shown
that grassmanians and partial flag manifolds have a toric degeneration that
can be decribed by a certain quiver. This type of quivers can be generalised to
what we call {\em flag quivers}.

We show in this paper that toric Gorenstein singularities $X$
provided by such flag quivers $\kQ$ are smoothable in codimension three,cf.\ Corollary \ref{cor-smoothable}.
\vspace{1ex}\\
The idea is to determine their infinitesimal deformation spaces $T^\kn_X$ 
($\kn=1,2$): 
$T_X^1$ describes the infinitesimal deformations,
and $T_X^2$ contains the obstructions for extending
deformations to larger base spaces -- see \cite{Loday} for more details.
The results will show that all deformations are unobstructed
(cf.\ Theorem \ref{th-lackobst})
and, moreover, that there are enough of them for providing a smoothing
in codimension three (cf.\ Theorem \ref{th-omega}).
\vspace{1ex}\\
The singularities $X$ turn out to be cones over singular
Fano varieties $\PP_{\kDQ}$. In Theorem \ref{th-multipath},
we determine the (embedded) infinitesimal deformations of these
projective varieties. 
This yields to a precise characterization of those flag quivers
$\kQ$ leading to Fanos $\PP_{\kDQ}$ which are smoothable in codimension three.
The condition is that every so-called simple knot
(cf.\ Definition \ref{def-simpknot}) can be by-passed with a multipath
connecting its neighbors,
cf.\ Corollary \ref{cor-smoothable} again.
\par

\neu{intro-contents}
The paper is organized as follows:
In \S \ref{tgor} we recall the main facts about
affine, toric Gorenstein singularities $X$
and their infinitesimal deformation
theory. In particular, we show how the modules $T^\kn_X$
are linked to certain invariants $\kT^\kn(\kDst)$ of the polytopes $\kDst$
defining $X$.
\vspace{1ex}
\\
In \S \ref{quiver},
we study mutually dual polytopes $\kDQ$ and $\kDstQ$ provided by quivers,
i.e.\ oriented graphs $\kQ$.
In the smooth case, i.e.\ if $\kDQ$ looks like an orthant in a
neighborhood
of every vertex, this has already been done in \cite{tilt} to
describe the moduli
spaces of representations of the quiver $\kQ$.
Here, we focus on the singular case.
\vspace{1ex}
\\
An important technical term is that of the tightness of a quiver $\kQ$.
This will be studied in \S \ref{tight}. Every quiver can be ``tightend'',
and this property straightens the relation between $\kQ$ ant its associated
polytopes. In particular, it allows 
to determine all
faces of $\kDstQ$ with a non-trivial $\kT^1$-invariant.
Then, we restrict ourselves to the special case of so-called flag quivers.
They are introduced (and their name is explained) in \S \ref{flag}. Moreover,
we determine their Picard group.
\vspace{1ex}
\\
Eventually, \S \ref{deform} deals with the deformation theory of flag quivers.
For them, it is possible to apply Theorem~\ref{th-TD} and
Theorem~\ref{th-D2van}, i.e.\ we can detect all non-vanishing $\kT^1$,
but prove the vanishing of $\kT^2$ for all relevant faces of $\kDstQ$.
This translates into the lack of obstructions for our singularities $X$.
\par

\section{Toric Gorenstein singularities}\label{tgor}

\neu{tgor-toric}
Let $N,M$ be two mutually dual free abelian groups of finite rank;
denote by $N_\R, M_\R$ the vector spaces obtained by extending the scalars.
Each polyhedral, rational cone $\sigma\subseteq N_\R$ with apex in $0$ 
gives rise to an affine, 
so-called toric variety $\toric{\sigma,N}:=\Spec\,\C[\sigma^\veee\cap M]$. 
See \cite{Da} for more details.\\
The toric variety $\toric{\sigma,N}$ is Gorenstein if and only if $\sigma$ is
the cone over a lattice polytope $\kP$ sitting in an affine hyperplane
of height one in $N_\R$,
i.e.\
if there is a primitive $R^\ast\in M$ such that 
$\Delta\subseteq [R^\ast=1]\subseteq N_\R$.
In this situation, we denote
$X_\kP:=\toric{\sigma,N}$.
The ring 
$A=\C\,[\sigma^\veee\cap M]$ as well as the
modules $T^\kn_{X_\kP}$ are $M$-graded, 
and the homogeneous pieces $T^\kn_{X_\Delta}(-R)$ with
$R\in M$ may be described in terms of the polytope $\kP$:
Consider the complex
\[
\xymatrix@!0@C=3.2em@R=7ex{
0 \ar[r] & C^0 \ar[rr] \ar@{=}[d] && C^1  \ar[rrr] \ar@{=}[d] &&&
C^2  \ar[rr] \ar@{=}[d] && C^3  \ar[r] \ar@{=}[d] & \dots\\
0 \ar[r] & N_\C  \ar[rr] && 
\makebox[7em]{$\oplus_{f_0\in\kP}
\raisebox{0.5ex}{$N_\C$}\hspace{-0.2em}\big/\hspace{-0.1em}
\raisebox{-0.5ex}{$\C\cdot f_0$}$}
\ar[rrr]&&& 
\makebox[7em]{$\oplus_{f_1<\kP} 
\raisebox{0.5ex}{$N_\C$}\hspace{-0.2em}\big/\hspace{-0.1em}
\raisebox{-0.5ex}{$\spann\,f_1$}$}
\ar[rr]&&\dots
}
\]
with $f_\kn<\kP$ in the definition of $C^{\kn+1}$ running through the 
$\kn$-dimensional faces of $\kP$; its cohomology is denoted by
$\kT^\kn(\kP):=H^\kn(C^\kb)$.
Then, in \cite{hodge} we have shown that
\par

\begin{theorem}[\cite{hodge}, (6.6)]\klabel{th-TD}
Assume that the two-dimensional faces of $\kP$ are either squares or 
triangles with area $1$ and $1/2$, respectively, i.e.\
$X_\kP$ is a conifold in codimension three. Then, if $R\in M$ is any
degree, we have for $\kn\leq 2$
\[
\renewcommand{\arraystretch}{1.2}
T^\kn_{X_\kP}(-R)=\left\{\begin{array}{@{}cl}
\kT^\kn\big(\kP\cap [R=1]\big) & \mbox{if } R\leq 1 \mbox{ on }\kP\\
0 & \mbox{otherwise}\,.
\end{array}\right.
\]
Moreover, the multiplication $\,x^s:T^\kn_X(-R)\to T^\kn_X(-(R-s))$
with $s\in\sigma^\veee\cap M$
is induced from the restriction of $\kT^\kn\big(\kP\cap [R=1]\big)$
to the face 
$\kP\cap [R-s=1]=
\kP\cap [R=1]\cap s^\bot$
of $\kP\cap [R=1]$.
\vspace{-2ex}
\end{theorem}
\par

\neu{tgor-van}
The vector space $\kT^1(\kP)$ or, to be more exact, a full-dimensional
polyhedral cone in it, parametrizes the set of Minkowski summands of
the polytope $\kP$:
Each vertex of a Minkowski summand is considered a mutation of an original
vertex $f_0\in\kP$ -- hence it provides an element of the
corresponding summand $N_\C/\C\cdot f_0$ of $C^1$.
While it is possible to determine $\kT^1(\kP)$, which leads to $T_X^1$, in many
cases,
we have to use a vanishing theorem for the
$\kT^2$-invariant which is responsible for the obstructions. 
In \cite{hodge} we have proved the following result:
\par

\ktrash{
{\bf Definition:}
We define an inductive process of {\em``cleaning'' vertices and 
two-dimensional faces} of $\kP$. At the beginning, all faces are 
assumed to be ``contaminated'', but then one may repeatedly apply the 
following rules (i) and (ii) in an arbitrary order:
\kitem
\item[(i)]
A two-dimensional $\an$-gon $\kPV{}<\kP$ is said to be clean if at least
$(\an-3)$ of its vertices are so. (In particular, every triangle is
automatically clean.)
\item[(ii)]
A vertex of $\kP$ is declaired to be clean if it is contained in no more than
$(n-3)$ two-dimensional faces that are not cleaned yet.
\kenditem
{\bf Theorem:}
{\em
Let $\kP$ be an $n$-dimensional, compact, convex polytope such that 
every three-dimensional face is a pyramid.
If every vertex 
(or, equivalently, every two-dimensional face)
may be cleaned in the sense of the previous definition, 
then $\kT^2(\kP)=0$.
}
\par
}

\begin{theorem}[\cite{hodge}, (1.1) and (4.7)]\klabel{th-D2van}
Let $\kP$ be an $n$-dimensional, compact, convex polytope such that every
three-dimensional face is a pyramid (with arbitrary base).
If no vertex is contained in more than $(n-3)$ two-dimensional,
non-triangular faces, then $\kT^2(\kP)=0$.
\vspace{-2ex}
\end{theorem}
\par

\section{Quiver polytopes}\klabel{quiver}

\neu{quiver-setup}
Let $\kQ$ be a connected quiver, i.e.\ an oriented graph.
It consists of a set $\kQ_0$ of vertices, a set $\kQ_1$ of arrows, and two
functions $\tail,\head:\kQ_1\to\kQ_0$ assigning to each arrow
$\alpha\in\kQ_1$ its tail $\tail(\alpha)$ and head $\head(\alpha)$.
This gives rise to the incidence matrix $\kA$; it consists of $\#\kQ_0$ rows
and $\#\kQ_1$ columns, and for $i\in \kQ_0$, $\alpha\in\kQ_1$ we have
\[
\kA_{i\alpha}:=\left\{
  \begin{array}{cl}
    +1&\mbox{if $i = \tail(\alpha)$ }\\
    -1&\mbox{if $i = \head(\alpha)$ }\\
    0&\mbox{otherwise}\,.
\end{array}
\right.
\]
The associated linear map 
$\pi: \Z^{\kQ_1}\to\Z^{\kQ_0}$, 
$\pi: [\alpha]\mapsto [{\tail(\alpha)}] - [{\head(\alpha)}]$ 
provides an exact sequence
\vspace{-1ex}
\[
0 \to \kF \stackrel{i}{\longrightarrow} \Z^{\kQ_1} 
\stackrel{\pi}{\longrightarrow} \Z^{\kQ_0} 
\stackrel{\ku{1}}{\longrightarrow} \Z 
\to 0
\]
with some free abelian group $\kF$ of rank $\#\kQ_1-\#\kQ_0+1$. It is
generated by the minimal, not necessarily oriented cycles in $\kQ$.
Denote by $\kH:=\ker(\Z^{\kQ_0}\to\Z)$ the so-called set of integral
weights; it contains a canonical one defined as 
$\wtc:=\pi(\ku{1})=\sum_{\alpha\in\kQ_1} \pi([\alpha])$.
\vspace{1ex}
\par

\begin{definition}\klabel{def-flow}
To any weight $\wt\in\kH_\R:=\kH\otimes_\Z\R$ we associate the so-called 
{\em flow polytope}
\vspace{-1ex}
\[
\kDQwt:= \pi^{-1}(\wt)\cap \R^{\kQ_1}_{\geq 0}\,.
\]
For non-connected quivers, the flow polytope is defined as the product of
the flow polytopes associated to the connected components.
\vspace{1ex}
\end{definition}

\begin{proposition}\klabel{prop-lattice}
If $\wt\in\kH$ is an integral weight, then $\kDQwt$ is a lattice polytope.
\vspace{-1ex}
\end{proposition}

\begin{proof}(Communicated by G.~M.~Ziegler)
The vertices of $\kDQwt$ may be obtained as the unique solutions of 
linear equations with a submatrix of $\kA$ as coefficients and
integers $\wt_p$ at the right hand side. On the other hand, square matrices
having
$\pm 1$ as the only non-trivial entries with each of them occuring at
most once in each column can only have determinant $\pm 1$ or $0$:\\
If every column contains both $1$ and $-1$, then the rows add up to zero.
Otherwise, there is a column having only one single non-trivial 
entry $\pm 1$ -- and we use exactly this one to develop our determinant
and end up with a smaller matrix.
\vspace{-1ex}
\end{proof}

\neu{quiver-weights}
There are two weights being of special interest:
\kitem
\item[(i)]
The {\em canonical weight} $\wtc:=\pi(\ku{1})$; we set $\kDQ:=\kDQwtc$. 
The shifted polytope
$\kDQ-\ku{1}= \kF_\R \cap \{ r\in\R^{\kQ_1}\,|\; r_\alpha\geq -1\}$ 
contains $0$ as an
interior lattice point. This makes it possible to define the dual
polytope as
\[
\kDQ^\veee:=\big\{a\in\kF_\R^\ast\,\big|\; \big\langle a,\,\kDQ-\ku{1}
\big\rangle\geq -1\big\}\,.
\] 
$\kDQ^\veee\subseteq \kF_\R^\ast$ is the smallest polytope 
containing $0\in\kF^\ast$ and the so-called {\em quiver polytope}
\[
\kDstQ:=\conv\big\{a^\alpha:=i^\ast([\alpha])\in\kF^\ast
\,\big|\;\alpha\in\kQ_1\big\}\,.
\]
In particular, both $\kDQ^\veee$ and $\kDstQ$ are compact lattice polytopes. 
\vspace{1ex}
\item[(ii)]
The {\em zero weight} $\wt:=0$. Then, $\kD(\kQ,0)\subseteq\kF$ 
is a polyhedral cone with apex; 
it is the dual cone of $\R_{\geq 0}\cdot\kDQ^\veee
=\R_{\geq 0}\cdot \kDstQ$.
Moreover, $\kD(\kQ,0)$ occurs as the 
``tail cone'' (cone of unbounded directions)
$\kD(\kQ,0)=\kDQwt^\infty$ for every weight $\wt\in\kH_\R$.
\kenditem
The quiver $\kQ$ lacks oriented cycles if and only if
\[
\kD(\kQ,0)=0 
\hspace{0.5em}\Longleftrightarrow\hspace{0.5em}
\kDQ \mbox{ is compact}
\hspace{0.5em}\Longleftrightarrow\hspace{0.5em}
0\in\innt\, \kDQ^\veee
\hspace{0.5em}\Longleftrightarrow\hspace{0.5em}
0\in\innt\, \kDstQ.
\]
If this is the case, then
$\kDstQ=\kDQ^\veee$ is a reflexive polytope in the sense of
Batyrev, cf.\ \cite{reflexive}. The corresponding affine
toric Gorenstein singularity
$X_\kDstQ:=\toric{\sigma,N}$ with 
$N:=\kF^\ast\oplus\Z$ and
$\sigma:=\cone\big(\kDstQ\big)\subseteq N_\R$
will be our main subject of investigation; it equals 
the cone over the projective toric variety $\PP_\kDQ:=
\toric{\normal(\kDQ), \kF^\ast}$
with $\normal(\kDQ)$ denoting the normal fan of $\kDQ$.
\vspace{2ex}

\begin{example}\klabel{ex-oct}
Let $\kQ$ be the quiver
\[
\unitlength=0.6pt
\begin{picture}(235.00,10.00)(28.00,790.00)
\put(43.00,800.00){\vector(1,0){60.00}}
\put(43.00,793.00){\vector(1,0){60.00}}
\put(123.00,800.00){\vector(1,0){60.00}}
\put(123.00,793.00){\vector(1,0){60.00}}
\put(203.00,800.00){\vector(1,0){60.00}}
\put(203.00,793.00){\vector(1,0){60.00}}
\put(33.00,796.50){\circle*{8}}
\put(113.00,796.50){\circle*{8}}
\put(193.00,796.50){\circle*{8}}
\put(273.00,796.50){\circle*{8}}
\end{picture}
\]
with the corresponding polytope  
\[
\kDstQ=\conv\{a^1,\dots,a^6\}\subseteq
\raisebox{1ex}{$\R^6$}\!\Big/
\raisebox{-1ex}{$\langle a^1+a^2,\;a^3+a^4,\;a^5+a^6\rangle$}
\]
being an octahedron with unit triangles as facets.
In particular, since $\kT^1$ and $\kT^2$ of an octahedron is $0$ and $\C^2$,
respectively,
Theorem \ref{th-TD} implies that $X_\kDstQ$ is rigid,
but $T^2_X=T^2(-R^\ast)$ is two-dimensional.
The singularity $X_\kDstQ$ is well known; its equations are the
$2$-minors of a general $(2\times 3)$-determinant. It illustrates the 
disappointing fact
that $T^2_X$ might contain more than just the obstructions,
cf.\ Remark~\zitat{deform}{T2}.
\end{example}

\neu{quiver-stability}
By a subquiver $\kQP\subseteq\kQ$ we mean a quiver $\kQP$ with $\kQP_0=\kQ_0$
and $\kQP_1\subseteq\kQ_1$. It is not assumed to be connected; in
particular there might even occur isolated points $p\in\kQP_0$.
Fixing a weight $\wt\in\kH_\R$, we will use the abbreviation 
$\wt(S):=\sum_{p\in S}\wt_p$ for subsets $S\subseteq\kQ_0$.
\par

\begin{definition}\klabel{def-stable}
Let $\wt\in\kH_\R$ be a weight of $\kQ$.
A subquiver $\kQP\subseteq\kQ$ is said to be
\kitem
\item
{\em $\wt$-(semi-) stable} (cf.\ \cite{King})
if any non-trivial, proper subset $S\subset\kQP_0$ that is
closed under $\kQP$-successors fulfills $\wt(S)<0$ (or $\leq 0$,
respectively);
\item
{\em $\wt$-polystable} if the connected components $\kQP^v$ of $\kQP$, meant
as connected quivers with a possibly smaller set $\kQP^v_0\subseteq\kQ_0$, 
fulfill $\wt(\kQP^v_0)=0$ and, moreover, are $\wt$-stable.
\kenditem
\end{definition}

While these notions where defined in \cite{King} to obtain decent moduli spaces
of quiver representations, we will just use them to describe the faces of our
quiver polytopes.
\vspace{1ex}

\begin{lemma}\klabel{lem-stable}
\begin{enumerate}
\item[(1)]
``Stable'' $\Longrightarrow$ ``polystable''
$\Longrightarrow$ ``semistable''.
\vspace{0.5ex}
\item[(2)]
Stable quivers are always connected. Semistability of 
$\kQP$ implies $\wt(\kQP^v_0)=0$ for its connected components.
\vspace{0.5ex}
\item[(3)]
The notions ``stability'' and ``semistability'' are closed under enlargement
of the subquiver $\kQP\subseteq\kQ$; ``polystability'' is not.
\vspace{0.5ex}
\item[(4)]
A subquiver $\kQP\subseteq\kQ$ is polystable if and only if $\kD(P,\wt)$
contains a point with positive coordinates, i.e.\ if the set
$\pi^{-1}(\wt)\cap \R^{\kQP_1}_{>0}$ is non-empty.
\vspace{0.5ex}
\item[(5)]
A subquiver $\kQP\subseteq\kQ$ is semistable if and only if 
$\kD(P,\wt)\neq\emptyset$.
\vspace{0.5ex}
\item[(6)]
Every semistable subquiver $\kQP\subseteq\kQ$ contains a 
(unique) maximal polystable subquiver $\bar{\kQP}\subseteq\kQP$.
\end{enumerate}
\end{lemma}

\begin{proof}
The first two parts are straightforward. Claim (3) follows from the easy fact
that the larger $\kQP\subseteq\kQ$, the more difficult is it for
$S\subset\kQ_0$ to be closed under $\kQP$-successors. However, the
corresponding
property fails for ``polystability'', since any enlargement of
$\kQP\subseteq\kQ$ may unify connected components.\\
To see that the condition in (4) suffices for polystability,
let $S\subseteq\kQ_0$ be an arbitrary subset.
We may use any
$r\in\pi^{-1}(\wt)$ to calculate $\wt(S)$ as
\[
\wt(S)= \sum_{S\stackrel{\alpha}{\rightarrow}(\kQ_0\setminus S)}
\hspace{-1em}r_\alpha -
\sum_{(\kQ_0\setminus S)\stackrel{\alpha}{\rightarrow} S}
\hspace{-1em}r_\alpha\,.
\]
Now, if $r\in\pi^{-1}(\wt)\cap \R^{\kQP_1}$ and $S$ is closed under 
$\kQP$-successors, then the first sum is void. However, if $S$ is not a
union of connected components of $\kQP$, then there must be at least
one $\kQP$-arrow
connecting $\kQ_0\setminus S$ and $S$, i.e.\
contributing to the second sum. In particular, if $r$ has only
positive coordinates, then $\wt(S)<0$.\\
For proving the necessity of the condition, we may assume that
$\kQP=\kQ$ is $\wt$-stable.
If
$\pi^{-1}(\wt)\cap \R^{\kQ_1}_{>0}=\emptyset$, then the vector
$\wt=(\wt_p)_{p\in\kQ_0}$ may not be written as a positive linear
combination of the columns of the incidence matrix $\kA$ introduced in
\zitat{quiver}{setup}. Thus, duality in convex geometry provides the existence
of
a non-trivial vector $h\in\R^{\kQ_0}/(\R\cdot\ku{1})$ 
such that 
$\langle h, \kA_{(\bullet,\alpha)}\rangle\geq 0$ for every arrow 
$\alpha\in\kQ_1$, but $\langle h,\wt\rangle\leq 0$.
The first property means $h_{\tail(\alpha)}\geq h_{\head(\alpha)}$. Hence,
denoting by $c_1<\dots<c_k$ ($k\geq 2$) the values of $h$ on $\kQ_0$ and
choosing an arbitrary $c_0<c_1$, the subsets
\[
S_v:=\{p\in\kQ_0\,|\; h_p \leq c_v\}\subset \kQ_0 \quad (v=0,\dots,k)
\]
are closed under successors. In particular, by the stability 
of $\kQ$, we obtain 
$\wt(S_v)<0$ for $v=1,\dots,k-1$ and
$\wt(S_0)=\wt(S_k)=0$.
This yields a contradiction via
\[
0 \,<\,
\sum_{v=1}^{k-1} (c_{v}-c_{v+1})\, \wt(S_v)
\,=\,\sum_{v=1}^k \;c_v\cdot \wt( S_v\setminus S_{v-1})
\,=\,\sum_{p\in\kQ_0}h_p\,\wt_p
\,\leq 0\,.
\]
Finally, we have to show (5) and (6). 
First, if $\kQP\subseteq\kQ$ is a
subquiver with $\kD(P,\wt)\neq\emptyset$, then we may choose an element
$r\in\kD(P,\wt)$ with maximal support 
$\bar{\kQP}_1:=\{\alpha\in\kQP_1\,|\; r_\alpha>0\}$. 
Hence, by (4), the corresponding subquiver $\bar{\kQP}$ 
is polystable and, using (1) and (3), $\kQP$ must be semistable.\\
It remains to check that $\kD(P,\wt)\neq\emptyset$ for semistable
$\kQP\subseteq\kQ$. We do this by copying the second part of the
proof of (4) with minor changes: The stronger condition
$\pi^{-1}(\wt)\cap \R^{\kQ_1}_{\geq0}=\emptyset$
implies the existence of an $h$ satisfying the strict inequality
$\langle h, \wt\rangle<0$. On the other hand, if $\kQ$ is semistable,
we have only $s_v\geq 0$. Nevertheless, one obtains the same contradiction.  
\end{proof}

\begin{corollary}\klabel{cor-polface}
For a subquiver $\kQP\subseteq\kQ$, we realize its flow polytope
as the subset
\[
\kDPwt
=\big\{r\in \kDQwt\,\big|\; r_\alpha=0 \,\mbox{ if }
\alpha\notin\kQP\big\}.
\]
This provides an order preserving
bijection between the poset of $\wt$-polystable subquivers, on the one hand,
and the face lattice of $\kDQwt$, on the other.
In particular, since the dimension of $\kDPwt$ is
$\big(\#\kQP_1 -\#\kQ_0 +\#(\mbox{\rm components of }\kQP)\big)$,
the $\wt$-polystable trees yield the vertices of $\kDQwt$.
\end{corollary}

\begin{proof}
Every face of $\kDQwt$ is of the form
$\kPf=\{r\in \kDQwt\,|\; r_{\alpha^i}=0\,,\, i=1,\dots,k\}$ for some edges
$\alpha^1,\dots,\alpha^k\in\kQ_1$. Assuming that the set
$\{\alpha^1,\dots,\alpha^k\}$ is maximal for the given face $\kPf$,
we obtain $\kQP$ by $\kQP_1:=\kQ_1\setminus \{\alpha^1,\dots,\alpha^k\}$.  
\end{proof}

\neu{quiver-faces}
Every connected quiver $\kQ$ is 
stable with respect to its {\em canonical weight} $\wtc$.
In this situation, we had defined in \zitat{quiver}{weights}
the polytopes $\kDstQ\subseteq\kDQ^\vee$. 
In general, under the dualization 
$\kD^\vee:=\{a \kst \langle a, \kD\rangle \geq -1 \}$
of polytopes containing the origin, we obtain an 
anti-isomorphism of the posets
\[
\xymatrix@!0@C=13em@R=3ex{
\Phi:
\{\mbox{faces of $\kD$ without $0$}\}
\ar[r]^-{\sim} &
\{\mbox{faces of $\kD^\vee$ without $0$}\}\\
\kPf \ar@{|->}[r] & \{a\in \kD^\vee \kst \langle a, \kPf\rangle = -1 \}.
}
\]
Note that faces containing $0$ correspond to faces of the dual tail cone
-- this gives a kind of a continuation of $\Phi$. 
If, e.g., $\kPf\leq\kD$ is as above and 
$[0,\kPf]\leq\kD$ denotes the smallest face
containing $0$ and $\kPf$, then
$\Phi([0,\kPf])$ is the tail cone of $\Phi(\kPf)\leq \kD^\vee$.
\vspace{0.5ex}\\
Applying this to $\kD:=\kDQ-\ku{1}$, we obtain for every $\wtc$-polystable
subquiver $\kQP\subseteq\kQ$
the face
\[
\kPF\big(\kQP,\kDstQ\big):= \Phi(\kDPwtc) = 
\conv \big\{a^\alpha\in \kF^\ast\,\big|\; \alpha\notin\kQP\big\}.
\]
Since $\kPF\big(\kQP,\kDstQ\big)$ does not contain $0$, it is also
a face of the quiver polytope $\kDstQ$.
\par

\begin{example}\klabel{ex-easy}
With $\kQ$ being the quiver\\
\hspace*{\fill}
\unitlength=1.2pt
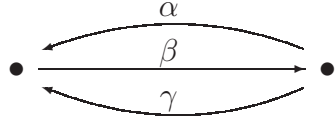
\begin{picture}(123.00,30.00)(120.00,770.00)
\put(99.00,788.00){\vector(1,0){84.00}}
\put(104.00,795.30){\vector(-3,-1){4.00}}
\put(104.00,781.00){\vector(-3,1){4.00}}
\put(92.00,788.00){\circle*{4}}
\put(190.00,788.00){\circle*{4}}
\put(140.00,807.00){\makebox(0,0)[cc]{$\alpha$}}
\put(140.00,793.00){\makebox(0,0)[cc]{$\beta$}}
\put(140.00,778.00){\makebox(0,0)[cc]{$\gamma$}}
\bezier300(101.00,794.50)(142.00,810.00)(182.00,794.50)
\bezier300(101.00,782.00)(144.00,765.00)(182.00,782.00)
\end{picture}
\\
we obtain
$\kF_\R^\ast=\R^3/\langle a-b+c\rangle$. Using the basis 
$\{a,c\}$, we can draw the polyhedra
$\kDQ\subseteq \kF_\R$ and $\kDQ^\vee,\kDstQ\subseteq \kF_\R^\ast$
as follows:
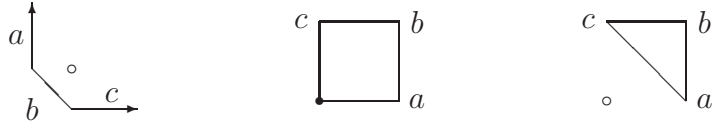
\begin{center}
\hspace*{\fill}
\unitlength=1.0pt
\begin{picture}(0,30)(0,0)
\put(0.00,15.00){\vector(0,1){25.00}}
\put(15.00,0.00){\vector(1,0){25.00}}
\put(0.00,15.00){\line(1,-1){15.00}}
\put(15.00,15.00){\circle{3}}
\put(30.00,5.00){\makebox(0,0)[cc]{$c$}}
\put(0.00,0.00){\makebox(0,0)[cc]{$b$}}
\put(-6.00,27.00){\makebox(0,0)[cc]{$a$}}
\end{picture}
\hspace*{\fill}
\unitlength=1.0pt
\begin{picture}(0,25)(0,-3)
\put(0.00,30.00){\line(1,0){30.00}}
\put(0.00,30.00){\line(0,-1){30.00}}
\put(30.00,0.00){\line(0,1){30.00}}
\put(0.00,00.00){\line(1,0){30.00}}
\put(0.00,0.00){\circle*{3}}
\put(37.00,0.00){\makebox(0,0)[cc]{$a$}}
\put(37.00,30.00){\makebox(0,0)[cc]{$b$}}
\put(-7.00,30.00){\makebox(0,0)[cc]{$c$}}
\end{picture}
\hspace*{\fill}
\unitlength=1.0pt
\begin{picture}(0,25)(0,-3)
\put(0.00,30.00){\line(1,0){30.00}}
\put(0.00,30.00){\line(1,-1){30.00}}
\put(30.00,0.00){\line(0,1){30.00}}
\put(0.00,0.00){\circle{3}}
\put(37.00,0.00){\makebox(0,0)[cc]{$a$}}
\put(37.00,30.00){\makebox(0,0)[cc]{$b$}}
\put(-7.00,30.00){\makebox(0,0)[cc]{$c$}}
\end{picture}
\hspace*{\fill}
\end{center}
Here are the five proper $\wtc$-polystable subquivers giving rise to faces
of them:
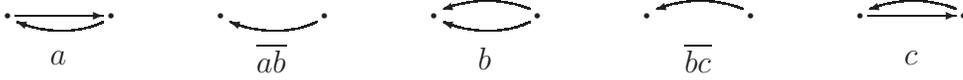
\begin{center}
\hspace*{\fill}
\unitlength=0.4pt
\begin{picture}(123.00,30.00)(120.00,770.00)
\put(99.00,788.00){\vector(1,0){84.00}}
\put(92.00,788.00){\circle*{4}}
\put(190.00,788.00){\circle*{4}}
\put(140.00,748.00){\makebox(0,0)[cc]{$a$}}
\put(104.00,781.00){\vector(-3,1){4.00}}
\bezier300(101.00,782.00)(144.00,765.00)(182.00,782.00)
\end{picture}
\hspace*{\fill}
\unitlength=0.4pt
\begin{picture}(123.00,30.00)(120.00,770.00)
\put(92.00,788.00){\circle*{4}}
\put(190.00,788.00){\circle*{4}}
\put(140.00,748.00){\makebox(0,0)[cc]{$\ko{ab}$}}
\put(104.00,781.00){\vector(-3,1){4.00}}
\bezier300(101.00,782.00)(144.00,765.00)(182.00,782.00)
\end{picture}
\hspace*{\fill}
\unitlength=0.4pt
\begin{picture}(123.00,30.00)(120.00,770.00)
\put(92.00,788.00){\circle*{4}}
\put(190.00,788.00){\circle*{4}}
\put(140.00,748.00){\makebox(0,0)[cc]{$b$}}
\put(104.00,795.30){\vector(-3,-1){4.00}}
\bezier300(101.00,794.50)(142.00,810.00)(182.00,794.50)
\put(104.00,781.00){\vector(-3,1){4.00}}
\bezier300(101.00,782.00)(144.00,765.00)(182.00,782.00)
\end{picture}
\hspace*{\fill}
\unitlength=0.4pt
\begin{picture}(123.00,30.00)(120.00,770.00)
\put(92.00,788.00){\circle*{4}}
\put(190.00,788.00){\circle*{4}}
\put(140.00,748.00){\makebox(0,0)[cc]{$\ko{bc}$}}
\put(104.00,795.30){\vector(-3,-1){4.00}}
\bezier300(101.00,794.50)(142.00,810.00)(182.00,794.50)
\end{picture}
\hspace*{\fill}
\unitlength=0.4pt
\begin{picture}(123.00,30.00)(120.00,770.00)
\put(99.00,788.00){\vector(1,0){84.00}}
\put(92.00,788.00){\circle*{4}}
\put(190.00,788.00){\circle*{4}}
\put(140.00,748.00){\makebox(0,0)[cc]{$c$}}
\put(104.00,795.30){\vector(-3,-1){4.00}}
\bezier300(101.00,794.50)(142.00,810.00)(182.00,794.50)
\end{picture}
\hspace*{\fill}
\end{center}

\end{example}

\neu{quiver-gamma}
Via contraction, we will construct new quivers $\kG_\kQ(\kQP)$
which allow to consider the faces 
$\kPF\big(\kQP,\kDstQ\big)\leq\kDstQ$ as autonomous quiver polytopes.

\begin{definition}\klabel{def-gamma}
For any subquiver $\kQP\subseteq\kQ$ we define a quiver
$\kG_\kQ(\kQP)$.
Its vertices $\kG_\kQ(\kQP)_0$ are the connected components
of $\kQP$, and the arrows are defined as
$\kG_\kQ(\kQP)_1:=\kQ_1\setminus\kQP_1$.
Every weight $\wt$ on $\kQ$ induces a weight $\wt$ on $\kG_\kQ(\kQP)$
with $\wtc$ staying the canonical weight.
If $\kQP$ was $\wt$-polystable, then $\wt$ turns into the $0$-weight
on $\kG_\kQ(\kQP)$.
\vspace{1ex}
\end{definition}

\begin{proposition}\klabel{prop-gamma}
Let $\kQP<\kQ$ be a non-empty, $\wtc$-polystable subquiver. Then, the face
$\kPF(\kQP,\kDstQ)<\kDstQ$ equals the quiver polytope
$\kDst\big(\kG_\kQ(\kQP)\big)$ and has dimension 
$\big(\#\kG_1-\#\kG_0\big)$.
Moreover, it is contained in a plane of $\kF^\ast$ having lattice
distance one from the origin.%
\vspace{-2ex}
\end{proposition}

\begin{proof}
Note that $\wtc=0$ in $\kG_\kQ(\kQP)$.
The original quiver $\kQ$ and
$\kG_\kQ(\kQP)$ are related by the following commutative diagram where
the vertical maps are surjective.
\[
\xymatrix@!0@C=5em@R=7ex{
0
\ar[r]
& \kF
\ar[r]
\ar[d]
& \Z^{\kQ_1}
\ar[r]
\ar[d]^-{\mbox{\rm\tiny forget}}
& \Z^{\kQ_0}
\ar[r]
\ar[d]^-{\mbox{\rm\tiny sum}}
& \Z
\ar[r]
\ar[d]
& 0\\
0
\ar[r]
& \kF(\kG)
\ar[r]
& \Z^{\kG_1}
\ar[r]
& \Z^{\kG_0}
\ar[r]
& \Z
\ar[r]
& 0
}
\]
Now, the first claim follows easily from the dual picture:
$\kF(\kG)^\ast\hookrightarrow\kF^\ast$ is saturated, and
$\kF(\kG)^\ast$ is the image of $\Z^{\kG_1}$ under the surjection
$\Z^{\kQ_1}\surj\kF^\ast$. 
The part concerning
the height is a consequence of the fact that the faces
of $\kD(\kQ)^\vee$ are contained in affine hyperplanes
$[r=-1]\subseteq\kF^\ast$ for certain vertices
$r\in \kD(\kQ)-\ku{1}$. By Proposition \ref{prop-lattice},
these $r$ are contained in the lattice $\kF$.
\vspace{-1ex}
\end{proof}

{\bf Example \ref{ex-easy}} (continued){\bf .}
The $\kDst(\kQ)$-faces corresponding to the five $\wtc$-polystable subquivers
are quiver polytopes arising from $\kG$ consisting of 
one vertex and one or two loops.

\section{Tightness}\klabel{tight}

\neu{tight-tight}
To ensure that there is a one-to-one
correspondence between arrows $\alpha\in\kQ_1$, on the one hand,
and  facets $\kD(\kQ\setminus\alpha,\,\wt)$ of $\kDQwt$,
on the other, we need the notion of tightness.
In particular, if $\kQ$ is $\wtc$-tight,
then $a^\alpha\in\kF^\ast$ will be the mutually distinct vertices of $\kDstQ$.

\begin{definition}\klabel{def-tight}
If $\wt\in\kH$ is an integral weight, then the quiver $\kQ$ is called
$\wt$-{\em tight} if for any $\alpha\in\kQ_1$ the subquiver
$\kQ\setminus\alpha$
is $\wt$-stable.
\vspace{2ex}
\end{definition}

\begin{lemma}\klabel{lem-tight}
\kitem
\item[(1)]
Let $\wt\in\kH$ 
such that $\kQ$ is $\wt$-stable.
By contraction of certain arrows in $\kQ$, the weight $\wt$ becomes
tight without changing the polytope $\kDQwt$.
Moreover, the canonical weight $\wtc$ may be tightened in such a way
that not only the polytope $\kDQwtc$, but also $\kDQwt$
remains untouched.
\vspace{0.5ex}
\item[(2)]
Assume that $\kQ$ is $\wt$-tight. If $\kQP\subseteq\kQ$ is a $\wt$-polystable
subquiver, then, $\kG_\kQ(\kQP)$ becomes $0$-tight.
\vspace{0.5ex}
\item[(3)]
Let $\kG$ be a $0$-tight quiver with $\#\kG_0\geq 2$.
Then, not counting the loops,
every knot of $\kG$ has at least two in- and outgoing arrows,
respectively. In particular, 
\[
\#\kG_1\;\geq\; 2\,\#\kG_0 + \#(\mbox{\rm loops of }\kG).
\vspace{-4ex}
\]
\kenditem
\end{lemma}

\begin{proof}
(1) If $\kQ\setminus\{\alpha\}$ is not $\wt$-stable, then there exists a
subset $S\subset\kQ_0$ that is, up to $\alpha$, closed under successors
and satisfies
$\wt(S)\geq 0$. Let $\beta_1,\dots,\beta_l$ be the arrows pointing from
$\kQ_0\setminus S$ to $S$; since $\kQ$ is $\wt$-stable, 
$\alpha$ leads in the opposite direction.
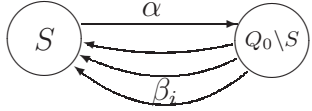
\begin{center}
\unitlength=0.7pt
\begin{picture}(163.00,45.00)(59.00,760.00)
\put(99.00,797.00){\vector(1,0){84.00}}
\put(104.00,786.00){\vector(-4,1){4.00}}
\put(106.00,776.50){\vector(-2,1){8.00}}
\put(102.00,769.00){\vector(-4,3){7.00}}
\put(79.00,789.00){\circle{40}}
\put(202.00,788.00){\circle{40}}
\put(79.00,789.00){\makebox(0,0)[cc]{$S$}}
\put(202.00,788.00){\makebox(0,0)[cc]{${\ks \kQ_0\setminus S}$}}
\put(137.00,805.00){\makebox(0,0)[cc]{$\alpha$}}
\put(144.00,760.00){\makebox(0,0)[cc]{$\beta_i$}}
\bezier80(102.00,786.00)(142.00,777.00)(180.00,786.00)
\bezier86(101.00,779.00)(144.00,759.00)(182.00,778.00)
\bezier113(97.00,772.00)(145.00,735.00)(188.00,771.00)
\end{picture}
\end{center}
If $r\in\pi^{-1}(\wt)$, then 
$r_\alpha=\sum_i r_{\beta_i}+\wt(S)\geq \sum_i r_{\beta_i}$. Hence, the
relations $r_{\beta_i}\geq 0$ together with the $\pi^{-1}(\wt)$-equations
force that $r_\alpha\geq 0$. In particular, we may contract $\alpha$
without changing the polytope $\kDQwt$ (including its lattice structure).
\vspace{0.5ex}\\
Tightening $\wtc$, we obtain $1-l=\wtc(S)\geq 0$. In case of $l=0$, the
$\wt$-stability of $\kQ$ implies that $\wt(S)>0$. In particular, 
contracting $\alpha$ does not change
neither $\kDQwtc$, nor $\kDQwt$. If $l=1$, then the situation is symmetric
with
$\wtc(S)=0$. Depending on whether $\wt(S)\geq 0$ or $\wt(S)\leq 0$, 
we should contract $\alpha$ or $\beta_1$, respectively.
\vspace{1.0ex}\\
(2) Let $\alpha\in \kG_1=\kQ_1\setminus\kP_1$.
Then, the connectivity of $\kQ\setminus\alpha$ implies
the connectivity of $\kG\setminus\alpha$.
Moreover, projecting any positive
$r\in \kD(\kQ\setminus\alpha,\theta)$ along the forgetful map
$\,\Z^{\kQ_1}\to\Z^{\kG_1}$ (see the diagram of the proof of 
Proposition \ref{prop-gamma}), provides
an $\ko{r}\in\kD(\kG\setminus\alpha,0)$ with positive entries.
\vspace{1.0ex}\\
(3) For every $\alpha\in \kG_1$, 
non-trivial subsets $S\subseteq\kG_0$ cannot be closed under
$(\kG\setminus\alpha)$-successors.
Hence, there is always a $\beta\in\kG_1\setminus\alpha$ leaving $S$.
Now, the claim follows from applying this fact to the cases
$\#S=1$ or $\#(\kG_0\setminus S)=1$.
\end{proof}

\neu{tight-smallfaces}
If $\kQ$ has no oriented cycles, then $\kDstQ=\kDQ^\vee$,
and Proposition \ref{prop-gamma} yields all its faces --
they equal $\kDst(\kG)$ with $\kG=\kG_\kQ(\kQP)$ for
$\wtc$-polystable subquivers $\kQP\leq\kQ$.
Moreover, $\kQ$ and hence $\kG$ maybe assumed to be $\wtc$-tight.
In particular,
$\#\kG_1\;\geq\; 2\,\#\kG_0 + \#(\mbox{\rm loops of }\kG)$.
\vspace{0.5ex}\\
Using this, we are now classifying 
all possible faces of those polytopes $\kDstQ$ 
up to dimension three. 
Note that $\wtc=0$ in $\kG$.
\par

{\em Dimension one:\,}
$\kG$ consists of one vertex with two loops. The corresponding polytope
$\kDst(\kG)$ is a lattice interval of length one.
\par

{\em Dimension two:\,}
The case $\#\kG_0=1$ yields the triple loop with $\kDst(\kG)$ being
the standard triangle. On the other hand, there is only one quiver 
with $\wtc(\kG)\equiv 0$
that consists of two vertices, four arrows, but no loops:
\begin{center}
\begin{picture}(150.00,20.00)(90.00,740.00)
\put(80.00,759.00){\vector(1,0){120.00}}
\put(80.00,753.00){\vector(1,0){120.00}}
\put(200.00,747.00){\vector(-1,0){120.00}}
\put(200.00,741.00){\vector(-1,0){120.00}}
\put(70.00,750.00){\circle{12}}
\put(212.00,750.00){\circle{12}}
\put(300.00,750.00){\makebox(0,0)[cc]{$\kGdbl(2)=\kGopp(2)$}}
\end{picture}
\end{center}
The corresponding polytope $\kDst(\kG)$ is the unit square.
\par

{\em Dimension three:\,}
The case $\#\kG_0\leq 2$ yields the quivers of the previous list with one
additional loop. Adding a loop to $\kG$ corresponds to taking the  
pyramid of height 1 over the corresponding polytope $\kDst(\kG)$. 
In particular, we obtain the unit
tetrahedron and the pyramid over the unit square.
\vspace{1ex}\\
On the other hand, there
are two different quivers involving three vertices and six edges:
\begin{center}
\hspace*{\fill}
\unitlength=0.8pt
\begin{picture}(110,90)(0,0)
\put(15.00,8.00){\vector(1,0){80.00}}
\put(15.00,12.00){\vector(1,0){80.00}}
\put(47.00,79.00){\vector(-2,-3){40.00}}
\put(50.00,76.50){\vector(-2,-3){39.00}}
\put(103.00,19.00){\vector(-2,3){40.00}}
\put(99.00,18.00){\vector(-2,3){39.00}}
\put(5.00,10.00){\circle{8}}
\put(105.00,10.00){\circle{8}}
\put(55.00,85.00){\circle{8}}
\put(120.00,70.00){\makebox(0,0)[cc]{$\kGdbl(3)$}}
\end{picture}
\hspace*{\fill}
\unitlength=0.8pt
\begin{picture}(110,90)(0,0)
\put(95.00,8.00){\vector(-1,0){80.00}}
\put(15.00,12.00){\vector(1,0){80.00}}
\put(7.00,19.00){\vector(2,3){40.00}}
\put(50.00,76.50){\vector(-2,-3){39.00}}
\put(63.00,79.00){\vector(2,-3){40.00}}
\put(99.00,18.00){\vector(-2,3){39.00}}
\put(5.00,10.00){\circle{8}}
\put(105.00,10.00){\circle{8}}
\put(55.00,85.00){\circle{8}}
\put(120.00,70.00){\makebox(0,0)[cc]{$\kGopp(3)$}}
\end{picture}
\hspace*{\fill}
\end{center}
The first quiver provides an octahedron. However, compared with the quiver
polytope presented in Example \zitat{quiver}{weights}, the central point does
no longer belong to the lattice. The other quiver provides the prism of height
$1$ over the unit triangle.
\par

\begin{corollary}\klabel{cor-two}
Let $\kQ$ be a quiver without oriented cycles. Then,
$\kDstQ$ and its faces satisfy the assumptions of
Theorem \ref{th-TD}: The two-dimensional faces are either squares or
triangles with area $1$ and $1/2$, respectively. Thus, 
$X_\kDstQ$ is a conifold in codimension three, and the vector spaces
$T^i_{X}$ may be obtained by calculating the corresponding
$\kT$-invariants of the $\kDstQ$-faces.
\end{corollary}

\neu{tight-T1}
Let $\kQ$ be a quiver without oriented cycles.
We conclude this chapter with determining all proper faces of $\kDstQ$
having a non-trivial $\kT^1$ (cf.\ \zitat{tgor}{toric}), 
i.e.\ admitting a non-trivial splitting into Minkowski summands.

\begin{example}\klabel{ex-perm}
If $\pi\in S_k$ is a permutation, then we denote by
$\kG(k,\pi)$ the quiver with
\kitem
\item
vertex set $\kG(k,\pi)_0=\Z/k\Z\,$ and
\vspace{0.5ex}
\item 
arrows $\beta_p,\gamma_p\in\kG(k,\pi)$ defined as $\beta_p:p\to (p+1)$ and 
$\gamma_p:p\to\pi(p)$ for $p=1,\dots,k$.
\vspace{-0.5ex}
\kenditem
The permutations $\pi^{\mbox{\tiny dbl}}(p):=p+1$ and 
$\pi^{\mbox{\tiny opp}}(p):=p-1$ are of special interest; the 
quivers $\kGdbl(k):=\kG(k,\pi^{\mbox{\tiny dbl}})$ and 
$\kGopp(k):=\kG(k,\pi^{\mbox{\tiny opp}})$ are 
double $n$-gons as shown in \zitat{tight}{smallfaces} for $k=2$ and $k=3$.
The corresponding polytopes are
\[
\kDst\big(\kGdbl(k)\big)=
[\mbox{crosspolytope of dimension $k$}]
\hspace{1em}\mbox{and}\hspace{1em}
\kDst\big(\kGopp(k)\big)=\Delta^{k-1}\times[0,1]\,.
\]
In particular, while $\kT^1\big(\kDst(\kGdbl(k))\big)=0$
for $k\geq 3$, we have
$\dim\kT^1\big(\kDst(\kGopp(k))\big)=1$ 
with the obvious Minkowski decomposition.
\end{example}

\begin{lemma}\klabel{lem-cycle}
Let $\kG$ be a quiver which is tight with respect to $\wtc=0$.
Assume that
$b\subseteq\kG_1$ is a simple cycle, i.e.\ not touching vertices twice.
\kitem
\item[(1)]
Contracting $b$, the resulting quiver $\bar{\kG}:=\kG/b$ is still tight with 
respect to
$\wtc= 0$. Moreover, $\kDst(\bar{\kG})$ is a face of 
$\kDst(\kG)$ inducing the restriction map
$p:\kT^1\big(\kDst(\kG)\big)\to\kT^1\big(\kDst(\bar{\kG})\big)$.
\vspace{0.5ex}
\item[(2)]
Unless $\kG=\kG(k,\pi)$ and $b$ is a cycle of length $k$, the map $p$ is
injective.
\kenditem
\end{lemma}

\begin{proof}
Let $b=\{\alpha^1,\dots,\alpha^k\}$ and denote by $a^i$ the image of
$[\alpha^i]$ in $\kF^\ast$.
The relations among the vertices
of $\kDst(\kG)$ which are induced from $\bar{\kG}$-knots are exactly the
relations among the vertices of $\kDst(\bar{\kG})$. On the other hand, the
knots being touched by $b$ express $(a^{i+1}-a^{i})$ as
an element of the vector space associated to the affine space 
$\A$ spanned by
$\kDst(\bar{\kG})$. In particular, $\kDst(\bar{\kG})$ is a face of
$\kDst(\kG)$ and, moreover, the remaining vertices $a^1,\dots,a^k$ are
contained in an affine plane being parallel to $\A$;
they form their own face $B:=\conv\{a^1,\dots,a^k\}$.
\vspace{1ex}\\
Let $\ku{t}\in\kT^1\big(\kDst(\kG)\big)$ -- here we think of it as a tuple
of dilatation factors for every compact edge of $\kDst(\kG)$:
The factors arise after applying the differential $d^1:C^1\to C^2$
from \zitat{tgor}{toric}; since it yields $0$, all the components of the image
must be contained in the subspaces $\spann f_1$.
\\
Since all vertical edges
connecting $\kDst(\bar{\kG})$ and $B$ have the same dilatation factor, we
may assume these factors to be zero. Now, if $p(\ku{t})=0$, then
the dilatation factors inside $\kDst(\bar{\kG})$ are also mutually equal;
it remains to show that they vanish.
If not, then we get a map
\[
\{a^1,\dots,a^k\}\longrightarrow \{\mbox{vertices of $\kDst(\bar{\kG})$}\}
\]
assigning $a^i$ the only vertex $a\in\kDst(\bar{\kG})$ such that
$\ko{a a^i}$ is an edge of $\kDst(\kG)$. Since this map is obviously surjective,
we obtain
\[
\#\bar{\kG}_1\leq k=\#b.
\]
Hence, $\kG$ equals $b$ with an additional arrow leaving and reaching each
knot.
\end{proof}

\begin{proposition}\klabel{prop-nontrivD1}
Let $\kQ$ be a $\wtc$-tight quiver without oriented cycles.
Then,
the only proper, $k$-dimensional faces $\kPF(\kQP,\kDstQ)$ of $\kDstQ$ 
having a non-trivial 
$\kT^1$ are those with $\kG_\kQ(\kQP)\cong\kGopp(k)$.  
\vspace{-1ex}
\end{proposition}

\begin{proof}
Let $\kG:=\kG_\kQ(\kQP)$. In the proof of the previous lemma 
we have seen that the existence
of a loop, i.e.\ of a cycle of length 1, implies that $\kDst(\kG)$ is a
pyramid over $\kDst(\bar{\kG})$. In particular, it has a trivial $\kT^1$.\\
On the other hand, via a successive contraction of simple cycles of 
length $k_i\geq 2$, we can 
produce a sequence of quivers, beginning with $\kG$, such that
\kitem
\item
we avoid the contraction of $k_i$-cycles in quivers isomorphic 
to $\kG(k_i,\pi)$, and
\vspace{0.5ex}
\item
we end with either the existence of loops or with a quiver isomorphic to
some $\kGdbl(m)$. The latter leads to a non-trivial $\kT^1$ only for $m=2$.
\vspace{-0.5ex}
\kenditem
By Lemma \ref{lem-cycle}(2), this sequence represents 
$\kT^1\big(\kDst(\kG)\big)$
as a subset of either $0$ or $\kT^1\big(\kDst(\kGdbl(2))\big)=\C$.
On the other hand, 
if the contraction of a simple $k_i$-cycle leads from $\kG^i$ to
${\kG}^{i+1}$, then
\[
\#({\kG}_0^{i+1})=\#(\kG_0^i) -(k_i-1)
\hspace{1em}\mbox{ and }\hspace{1em}
\#({\kG}_1^{i+1})=\#(\kG_1^i) -k_i\,.
\]
In particular, relations like $\#(\kG_1^i)\geq2\,\#(\kG_0^i)$ or
$\#(\kG_1^i)>2\,\#(\kG_0^i)$ 
will be inherited from $\kG^i$ to ${\kG}^{i+1}$. If $k_i\geq 3$, then
the weak inequality for $\kG^i$ does even imply the strict one for
${\kG}^{i+1}$. Thus, if our sequence ends with $\kGdbl(2)$,
only $2$-cycles have been contracted successively from $\kG$. 
This enforces that $\kG\cong\kGopp(k)$.
\end{proof}

\section{Flag Quivers}\klabel{flag}

\neu{flag-div}
First, we will
describe the classes of Weil and Cartier divisors on the
projective variety $\PP_\kDQwt$ provided by a general
quiver $\kQ$ with weight $\wt\in\kH$.
Assume, w.l.o.g., that $\kQ$ is $\wt$-tight.
We introduce the notation
\[
\kwtTrees(\kQ):=\{\wt\mbox{-polystable trees }\kQT<\kQ\}.
\]
The tightness of $\kQ$ implies that $\bigcap\kwtTrees(\kQ)=\emptyset$. 
For $\kQT\in\kwtTrees(\kQ)$, or more general for any $\wt$-polystable
subquiver $\kQP<\kQ$, we may define
\[
\kH(\kQP):= \kH\big(\kG_\kQ(\kQP)\big)
=\ker\Big( \Z^{\mbox{\tiny $\kQP$-comp}}
\stackrel{\mbox{\tiny deg}}{\longrightarrow}\Z\Big)
\]
with $\kH$ and $\kG_\kQ(\kQP)$ as mentioned in \zitat{quiver}{setup} and 
Definition \ref{def-gamma}, respectively.
\vspace{1ex}

\begin{proposition}\klabel{prop-Weil}
The class group $\mbox{\rm DivCl}\,(\PP_\kDQwt)$ of Weil divisors 
equals $\kH$.
\vspace{-1ex}
\end{proposition}

\begin{proof}
Equivariant Weil divisors correspond to maps $\CN^{(1)}\to\Z$ with
$\CN$ denoting the inner normal fan $\normal(\kDQwt)$ of
$\kDQwt\subseteq\kF_\R$. 
Since the elements of
$\CN^{(1)}$ encode the facets of $\kDQwt$, i.e.\ the $\wt$-polystable
subquivers $\kQ\setminus\alpha$, the equivariant Weil divisors 
may be written as elements of 
$\Z^{\kQ_1}$.\\
On the other hand, as being well known in the theory of toric varieties,
the principal divisors among the equivariant ones are given by
$\kF\subseteq\Z^{\kQ_1}$. Hence, the claim follows from
\[
\raisebox{0.5ex}{$\Z^{\kQ_1}$}\hspace{-0.3em}\big/
\hspace{-0.2em}\raisebox{-0.5ex}{$\kF$}
=
\im\Big(  \Z^{\kQ_1}\to \Z^{\kQ_0}\Big)
=
\ker\Big( \Z^{\kQ_0} \stackrel{\mbox{\tiny deg}}{\longrightarrow}\Z\Big)
=
\kH\,.
\vspace{-4ex}
\]
\end{proof}

\neu{flag-Cartier}
Using the fact that divisors on a toric variety are locally principal 
if and only if they are principal on the equivariant affine charts, we
obtain a description of $\Pic\,(\PP_\kDQwt)$ as well.
\par

\begin{proposition}\klabel{prop-Pic}
The Picard group $\Pic\,(\PP_\kDQwt)$ of the projective toric variety
$\PP_\kDQwt$ equals 
\[
\Pic\,(\PP_\kDQwt)=\ker\Big(
\kH \rightarrow
\oplus_{\kQT\in\kwtTrees}\kH(\kQT)\Big)\,.
\vspace{-1ex}
\]
\end{proposition}

\begin{proof}
An element $g\in\Z^{\kQ_1}$ represents a Cartier divisor if and only if
$g$ induces principal divisors on the affine charts $\toric{\delta(\kQT)}$
with $\delta(\kQT):=\langle a^\alpha\in\kF^\ast\,|\; \alpha\notin\kQT\rangle
\in\normal\big(\kDQwt\big)$ for every $\wt$-polystable tree $\kQT\in\kwtTrees$.
This condition is equivalent to $g\in\Z^{\kQT_1}+\kF\subseteq\Z^{\kQ_1}$.
Adapting the commutative diagram from the proof of
Proposition \ref{prop-gamma}, we obtain
\vspace{-2ex}
\[
\hspace*{5em}
\xymatrix@!0@C=5em@R=7ex{
&0\ar[d]&0\ar[d]&&\\
0\ar[r]&\Z^{\kQT_1}\cap\kF\ar[r]\ar[d]&\Z^{\kQT_1}\ar[d]&&\\
0
\ar[r]
& \kF
\ar[r]
\ar[d]
& \Z^{\kQ_1}
\ar[r]
\ar[d]^-{\mbox{\rm\tiny forget}}
& \kH
\ar[r]
\ar[d]^-{\mbox{\rm\tiny sum}}
& 0\\
0
\ar[r]
& \kF(\kG)
\ar[r]\ar[d]
& \Z^{\kG_1}
\ar[r]\ar[d]
& \kH(\kQT)
\ar[r]\ar[d]
& 0\\
&0&0&0
}
\]
It implies that the membership $g\in\Z^{\kQT_1}+\kF$ translates into
the fact that the class $\bar{g}\in\kH$ maps to $0\in\kH(\kQT)$.
\vspace{1ex}
\end{proof}

\begin{corollary}\klabel{cor-Pic}
$\Pic\,(\PP_\kDQwt)=
\begin{array}[t]{@{\,}r}
\big\{\wt^\prime\in\kH\kSt
\wt^\prime(S)=0 \mbox{\rm\ for }
S\subseteq\kQ_0 \mbox{\rm\ with }
\wt(S)=0  
\hspace{3em}\vspace{0.2ex}\\ 
\mbox{\rm and }\kQ_{|S} \mbox{\rm\ and } \kQ_{|(\kQ_0\setminus S)}
\mbox{\rm\ being $\wt$-semistable}\big\}.
\end{array}
$%
\vspace{0.8ex}
\\
In particular, a necessary, but only on $\wt\in\kH$
depending condition for $\Pic\,(\PP_\kDQwt)\cong\Z$
is that
$
\;\big\{\wt^\prime\in\kH_\R\kst
\wt^\prime(S)=0 \mbox{\rm\ for } S\subseteq\kQ_0 \mbox{\rm\ with }
\wt(S)=0\big\}=\R\cdot\wt\,.
$
\end{corollary}

\neu{flag-setup}
Now, we turn to the main point of this section -- the introduction of
the so-called flag quivers. They will allow an easy description of their
Picard group as well as, in \S \ref{deform}, of their deformation theory.

\begin{definition}\klabel{def-flagq} 
A quiver $\kQ$ without oriented cycles 
is called a {\em flag quiver} if
\kitem
\item[(i)]
there is exactly one source $p^0\in\kQ_0$ with $m:=\wtc(p^0)$,
\vspace{0.5ex}
\item[(ii)]
there are sinks $p^1,\dots,p^l$ with $m_i:=-\wtc(p^i)\geq 2$, and
\vspace{0.5ex}
\item[(iii)]
the canonical weight vanishes on the remaining knots, i.e.\
on $\kQ_0\setminus\{p^0,\dots,p^l\}$. In particular, we have
$m=\sum_{i=1}^lm_i$.
\vspace{1.0ex}
\kenditem
\end{definition}

\begin{remark}
The condition ``$m_i\geq 2$'' may be explained as follows: If
$m_i=1$,
then $\kQ$ cannot be tight, cf.\ Lemma \ref{lem-polytrees}.
Moreover, tightening would
mean to contract the arrow
pointing to $p^i$, hence to create a non-sink with negative weight.
\end{remark}

\neu{flag-bcks}
The name ``flag quiver'' arises from the following example from
\cite{flags}:\\
For positive integers $k_1,\dots,k_{l+1}$, we set 
$n_i:=\sum_{v=1}^i k_v$
and $n:=n_{l+1}$. 
Then, a quiver $\kQ(n_1,\dots,n_l,n)$ may be defined via an
$(n\times n)$-ladder-box containing the $(k_i\times k_i)$-squares on the main
diagonal. As depicted in the middle figure below,
$\kQ_0$ consists of the interior lattice points in the region below
the small squares and of the closures of the
$(l+1)$ connected components of the part of the
boundary of this region that avoids the south west corners of the
$(k_i\times k_i)$-squares. As arrows we take all possibilities pointing 
eastbound and northbound.
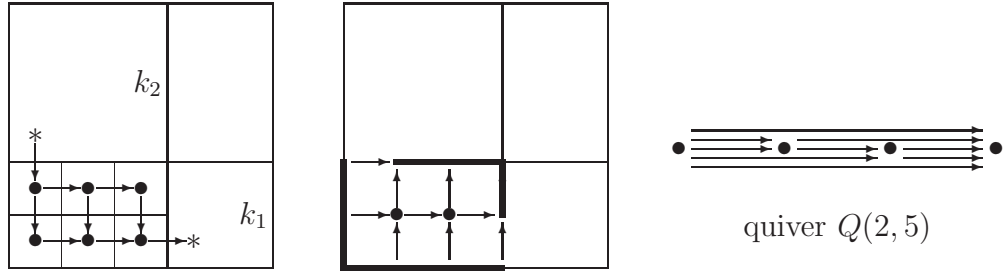
\begin{center}
\hspace*{\fill}
\unitlength=0.5pt
\begin{picture}(200.00,200.00)(40.00,640.00)
\put(40.00,840.00){\line(0,-1){200.00}}
\put(40.00,640.00){\line(1,0){200.00}}
\put(240.00,640.00){\line(0,1){200.00}}
\put(240.00,840.00){\line(-1,0){200.00}}
\put(40.00,720.00){\line(1,0){200.00}}
\put(160.00,840.00){\line(0,-1){200.00}}
\put(60.00,740.00){\makebox(0,0)[cc]{$\ast$}}
\put(180.00,660.00){\makebox(0,0)[cc]{$\ast$}}
\put(60.00,700.00){\makebox(0,0)[cc]{$\bullet$}}
\put(100.00,700.00){\makebox(0,0)[cc]{$\bullet$}}
\put(140.00,700.00){\makebox(0,0)[cc]{$\bullet$}}
\put(60.00,660.00){\makebox(0,0)[cc]{$\bullet$}}
\put(100.00,660.00){\makebox(0,0)[cc]{$\bullet$}}
\put(140.00,660.00){\makebox(0,0)[cc]{$\bullet$}}
\put(60.00,734.00){\vector(0,-1){28.00}}
\put(100.00,694.00){\vector(0,-1){28.00}}
\put(140.00,694.00){\vector(0,-1){28.00}}
\put(60.00,694.00){\vector(0,-1){28.00}}
\put(66.00,700.00){\vector(1,0){28.00}}
\put(106.00,700.00){\vector(1,0){28.00}}
\put(66.00,660.00){\vector(1,0){28.00}}
\put(106.00,660.00){\vector(1,0){28.00}}
\put(146.00,660.00){\vector(1,0){28.00}}
\put(225.00,680.00){\makebox(0,0)[cc]{$k_1$}}
\put(145.00,780.00){\makebox(0,0)[cc]{$k_2$}}
\linethickness{0.1pt}
\put(80.00,640.00){\line(0,1){80.00}}
\put(120.00,640.00){\line(0,1){80.00}}
\put(40.00,680.00){\line(1,0){120.00}}
\end{picture}
\hfill
\begin{picture}(200.00,200.00)(40.00,640.00)
\put(40.00,840.00){\line(0,-1){200.00}}
\put(40.00,640.00){\line(1,0){200.00}}
\put(240.00,640.00){\line(0,1){200.00}}
\put(240.00,840.00){\line(-1,0){200.00}}
\put(80.00,720.00){\line(1,0){160.00}}
\put(160.00,840.00){\line(0,-1){160.00}}
\put(80.00,680.00){\makebox(0,0)[cc]{$\bullet$}}
\put(120.00,680.00){\makebox(0,0)[cc]{$\bullet$}}
\put(46.00,720.00){\vector(1,0){28.00}}
\put(46.00,680.00){\vector(1,0){28.00}}
\put(86.00,680.00){\vector(1,0){28.00}}
\put(126.00,680.00){\vector(1,0){28.00}}
\put(80.00,686.00){\vector(0,1){28.00}}
\put(120.00,686.00){\vector(0,1){28.00}}
\put(160.00,686.00){\vector(0,1){28.00}}
\put(80.00,646.00){\vector(0,1){28.00}}
\put(120.00,646.00){\vector(0,1){28.00}}
\put(160.00,646.00){\vector(0,1){28.00}}
\linethickness{2pt}
\put(78.00,720.00){\line(1,0){84.00}}
\put(160.00,722.00){\line(0,-1){44.00}}
\put(40.00,722.00){\line(0,-1){84.00}}
\put(38.00,640.00){\line(1,0){124.00}}
\end{picture}
\hfill
\begin{picture}(240.00,160.00)(40.00,640.00)
\put(50.00,737.00){\vector(1,0){60.00}}
\put(50.00,730.00){\vector(1,0){60.00}}
\put(50.00,723.00){\vector(1,0){140.00}}
\put(130.00,737.00){\vector(1,0){140.00}}
\put(130.00,730.00){\vector(1,0){60.00}}
\put(210.00,730.00){\vector(1,0){60.00}}
\put(210.00,723.00){\vector(1,0){60.00}}
\put(50.00,744.00){\vector(1,0){220.00}}
\put(50.00,716.00){\vector(1,0){220.00}}
\put(40.00,730.00){\makebox(0,0)[cc]{$\bullet$}}
\put(120.00,730.00){\makebox(0,0)[cc]{$\bullet$}}
\put(200.00,730.00){\makebox(0,0)[cc]{$\bullet$}}
\put(280.00,730.00){\makebox(0,0)[cc]{$\bullet$}}
\put(160.00,670.00){\makebox(0,0)[cc]{quiver $\kQ(2,5)$}}
\end{picture}
\hspace*{\fill}
\end{center}
\par
In \cite{flags}, the authors have originally considered
$\kQ^\ast(n_1,\dots,n_l,n)$ as depicted in the left figure above:
It is a kind of a dual quiver;
its ordinary vertices correspond to the boxes instead 
of the lattice points, and there are additional, exeptional, vertices called
``$\ast$'' sitting in the south west corners of the
$(k_i\times k_i)$-squares. Nevertheless,
it was shown that the corresponding $X_\kDstQ$
equals the cone over a projective toric variety being a degeneration of the flag
manifold $\mbox{Flag}(n_1,\dots,n_l,n)$. The polytopes assigned to
the quiver $\kQ(n_1,\dots,n_l,n)$ are called
$\kD(n_1,\dots,n_l,n)$ and $\kDst(n_1,\dots,n_l,n)$, respectively.
\par

\neu{flag-poly}
The polystability of subquivers has an easy meaning in the special case
of flag quivers:

\begin{lemma}\klabel{lem-polytrees}
Let $\kQ$ be a flag quiver.
$\kQP\subseteq\kQ$ is $\wtc$-polystable iff 
it is a union of paths leading from $p^0$ to 
{\em every} sink $p^i$ ($i=1,\dots,l$).
Moreover, $\kQ$ is tight if and only if there are no verices with only
one in- and outgoing arrow, respectively. In particular, tightening preserves
the property of being a flag quiver.
\end{lemma}

\begin{proof}
Both the criterion for polystability and the neccessity of the tightness
condition for $\kQ$ are clear.\\
On the other hand, assume that $\kQ$ 
satisfies this condition and let $\alpha\in\kQ_1$ be an arbitrary arrow.
We may choose paths $r^v$ leading from $p^0$ to $p^v$, but avoiding
$\alpha$. Moreover, for any $\beta\neq\alpha$ there is a special path
$s(\beta)$ which, additionally, touches $\beta$; let $p^{\mu(\beta)}$ be the
sink reached by $s(\beta)$. Then, with
\[
R(\beta):= s(\beta) \cup \Big(\mbox{$\bigcup$}_{v\neq\mu(\beta)}r^v\Big)
\]
we have obtained a union of paths encoding a polystable subquiver containing
$\beta$, but not $\alpha$. In particular, 
$\bigcup_{\beta\neq\alpha} R(\beta)$ shows the polystability of the quiver
$\kQP=\kQ\setminus\alpha$.
\end{proof}

Tight flag quivers may be visualized as a so-called {\em \fence}, 
i.e.\ as a system
of mutually intersecting ropes leading from the $l$ different
ceilings $p^1,\dots,p^l$ to the only base $p^0$. The intersections correspond
to the knots $b\in\kQ_0\setminus\{p^0,\dots,p^l\}$.
If, moreover, the quiver is a plane one, then,
by Corollary \ref{cor-polface},
the dimension of the 
corresponding polytopes
$\kDQ$ and $\kDstQ$ may be read of the plane \fence as the number of 
compact regions.

\begin{example}\klabel{ex-fence}
Here, we present three fences of dimensions six, five, and again five.
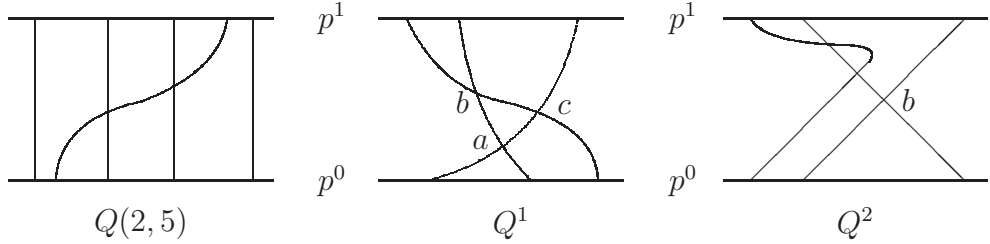
\begin{center}
\hspace*{\fill}
\unitlength=0.5pt
\begin{picture}(220.00,140.00)(40.00,660.00)
\thicklines
\put(40.00,800.00){\line(1,0){200.00}}
\put(40.00,677.00){\line(1,0){200.00}}
\thinlines
\bezier106(75.00,677.00)(79.00,724.00)(140.00,737.00)
\bezier106(140.00,737.00)(200.00,760.00)(205.00,800.00)
\put(225.00,800.00){\line(0,-1){123.00}}
\put(60.00,800.00){\line(0,-1){123.00}}
\put(115.00,800.00){\line(0,-1){123.00}}
\put(165.00,800.00){\line(0,-1){123.00}}
\put(140.00,645.00){\makebox(0,0)[cc]{$\kQ(2,5)$}}
\put(285.00,677.00){\makebox(0,0)[cc]{$p^0$}}
\put(285.00,800.00){\makebox(0,0)[cc]{$p^1$}}
\end{picture}
\hspace*{\fill}
\unitlength=0.5pt
\begin{picture}(200.00,140.00)(40.00,660.00)
\thicklines
\put(40.00,800.00){\line(1,0){185.00}}
\put(40.00,677.00){\line(1,0){185.00}}
\thinlines
\bezier106(60.00,800.00)(86.00,746.00)(133.00,737.00)
\bezier106(133.00,737.00)(204.00,720.00)(205.00,677.00)
\bezier106(100.00,800.00)(110.00,720.00)(155.00,677.00)
\bezier106(75.00,677.00)(170.00,700.00)(190.00,800.00)
\put(140.00,645.00){\makebox(0,0)[cc]{$\kQ^1$}}
\put(117.00,706.00){\makebox(0,0)[cc]{$a$}}
\put(180.00,733.00){\makebox(0,0)[cc]{$c$}}
\put(103.00,737.00){\makebox(0,0)[cc]{$b$}}
\put(270.00,677.00){\makebox(0,0)[cc]{$p^0$}}
\put(270.00,800.00){\makebox(0,0)[cc]{$p^1$}}
\end{picture}
\hspace*{\fill}
\unitlength=0.5pt
\begin{picture}(200.00,140.00)(40.00,660.00)
\thicklines
\put(40.00,800.00){\line(1,0){200.00}}
\put(40.00,677.00){\line(1,0){200.00}}
\thinlines
\put(100.00,677.00){\line(1,1){123}}
\put(60.00,677.00){\line(1,1){90}}
\bezier106(60.00,800.00)(70.00,782.00)(125.00,780.00)
\bezier106(125.00,780.00)(160.00,780.00)(150,767)
\put(100.00,800.00){\line(1,-1){123}}
\put(180.00,738.00){\makebox(0,0)[cc]{$b$}}
\put(140.00,645.00){\makebox(0,0)[cc]{$\kQ^2$}}
\end{picture}
\hspace*{\fill}
\vspace{1ex}
\end{center}
\end{example}

\neu{flag-onesource}
Assume that $\kQ$ is a tight flag quiver.
Denoting by $B\subseteq\kQ_0\setminus\{p^0,\dots,p^l\}$ the set of
{\em blocking knots}, i.e.\ those that are not avoidable in 
a set of paths leading from $p^0$ to each of the ends $p^1,\dots,p^l$,
respectively,
the Picard number of $\PP_\kDQ$ will be $\#(B)+l$. More precisely,
we  obtain

\begin{proposition}\klabel{prop-flagPic}
$
\;\Pic\,(\PP_\kDQ)=
\ker\big(\Z^{\{p^0,\dots,p^l\}\cup B}
\stackrel{\mbox{\tiny\rm deg}}{\longrightarrow}\Z\big)\,.
$
\vspace{-1ex}
\end{proposition}

\begin{proof}
The $\wtc$-polystable subquivers are the unions of paths leading from
$p^0$ to $\{p^1,\dots,p^l\}$. In particular, if
$\wt^\prime\in\Pic\,(\PP_\kDQ)$, then $\wt^\prime(b)=0$ for vertices
$b\notin\{p^0,\dots,p^l\}\cup B$. On the other hand, since for each
$S$ as in Corollary \zitat{flag}{Cartier} either $S$ or $\kQ_0\setminus S$
contains $\{p^0,\dots,p^l\}\cup B$, this remains the only condition.
\end{proof}

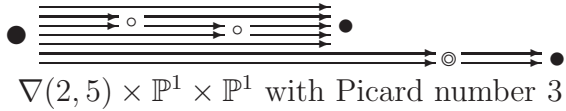
\begin{example}\klabel{ex-ones}
\begin{center}
\hspace*{\fill}
\unitlength=0.5pt
\begin{picture}(400.00,50.00)(40.00,720.00)
\put(50.00,737.00){\vector(1,0){60.00}}
\put(50.00,730.00){\vector(1,0){60.00}}
\put(50.00,723.00){\vector(1,0){140.00}}
\put(130.00,737.00){\vector(1,0){140.00}}
\put(130.00,730.00){\vector(1,0){60.00}}
\put(210.00,730.00){\vector(1,0){60.00}}
\put(210.00,723.00){\vector(1,0){60.00}}
\put(50.00,744.00){\vector(1,0){220.00}}
\put(50.00,716.00){\vector(1,0){220.00}}
\put(50.00,709.00){\vector(1,0){300.00}}
\put(50.00,702.00){\vector(1,0){300.00}}
\put(370.00,709.00){\vector(1,0){60.00}}
\put(370.00,702.00){\vector(1,0){60.00}}
\put(33.00,723.00){\circle*{14}}
\put(120.00,733.00){\circle{6}}
\put(200.00,727.00){\circle{6}}
\put(282.00,730.00){\circle*{10}}
\put(360.00,705.00){\circle{6}}
\put(360.00,705.00){\circle{9}}
\put(442.00,705.00){\circle*{10}}
\put(240.00,680.00){\makebox(0,0)[cc]{$\kD(2,5)\times\PP^1\times\PP^1$ with
Picard number 3}}
\end{picture}
\hspace*{\fill}
\vspace{4ex}
\end{center}
Knots $x\in B$ give rise
to a decomposition of $\kQ$ into a join of smaller quivers,
meaning that $\PP_\kDQ$ splits into a product of projective varieties.
\end{example}

\section{Deformation theory of flag quivers}\label{deform}

\neu{deform-T1}
Let $\kQ$ be a tight flag quiver.
From Theorem \ref{th-TD} and Corollary \ref{cor-two} we know that 
the module $T^1_{X}$ of infinitesimal deformations of
$X_\kDstQ=\Cone(\PP_\kDQ)$
is built from the spaces of Minkowski summands
$\kT^1(\kPF)$ of the faces $\kPF\leq\kDstQ$.
These faces look like
$\kPF(\kQP,\kDstQ)=\kDst(\kG_\kQ(\kQP))$
for $\wtc$-polystable subquivers $\kQP\leq\kQ$,
and, by Proposition \ref{prop-nontrivD1},
$\kG_\kQ(\kQP)$ must be isomorphic to $\kGopp(k)$
to yield a non-trivial $\kT^1$.
On the other hand, in the special case of flag quivers,
Lemma \ref{lem-polytrees} provides an explicit description
of the $\wtc$-polystable subquivers at all. Combining all this information,
we will get a complete description of $T^1_{X}$.

\begin{definition}\klabel{def-simpknot}
A knot $b\in\kQ_0\setminus\{p^0,\dots,p^l\}$ in a tight flag quiver
is called {\em simple} if 
\vspace{0.5ex} 
\kitem
\item
$b$ is
of valence four, i.e.\
supporting exactly two pairs of in- and outgoing arrows, 
\vspace{1ex}
respectively, and
\item
besides $b$,
both pairs do neither have 
a common tail or head of valence four in the set
$\kQ_0\setminus\{p^0,\dots,p^l\}$, 
neither a common head $p^i$ with $m_i=2$,
nor the common tail $p^0$ with $m=2$. 
\vspace{1ex}
\kenditem
Visualizing $\kQ$ as a \fence, then simple knots correspond to
the simple crossings of two ropes that are not adjacent to a further
simple crossing
of the same two ropes.
\end{definition}

{\bf Example \ref{ex-fence}} (continued){\bf .}
In the first two quivers $\kQ(2,5)$ and $\kQ^1$ of Example
\ref{ex-fence}, every knot of $\kQ\setminus \{p^0,p^1\}$ is simple.
On the other hand, in $\kQ^2$, only $b$ shares this property. The remaining
two knots provide for each other the reason to violate the condition
described in the previous definition.

\begin{proposition}\klabel{prop-simpknot}
Let $\kQ$ be a tight flag quiver with $\dim\kDstQ\geq 3$. 
Then, the only faces $\kPF\leq\kDstQ$ admitting
a non-trivial Minkowsi decomposition are 
the two-dimensional squares $\kPF(\kQ\setminus b)$ with
$b$ being a simple knot.\\
{\rm ($\kQ\setminus b\subseteq\kQ$ denotes the subquiver obtained
by erasing the four arrows containing $b$.)
}
\end{proposition}

\begin{proof}
First, we consider the {\em proper} faces of $\kDstQ$. 
Lemma \ref{lem-polytrees} 
tells us that $\wtc$-polystable subquivers $\kQP$
consist of one
major component and a bunch of isolated knots from
$\kQ_0\setminus\{p^0,\dots,p^l\}$. 
On the other hand, the only quivers providing
a non-trivial $\kT^1$ are $\kGopp(k)$. If $\kG_\kQ(\kQP)=\kGopp(k)$
with $k\geq 3$, then $\kQ$ would have to contain oriented cycles.
Thus, $\kQP=\kQ\setminus\{b\}$ for some knot $b$.
\vspace{1ex}\\
Denote by 
$\alpha^1,\alpha^2$ and $\beta^1,\beta^2$ the pairs of in- and outgoing
$b$-arrows, respectively. 
Assuming that, for instance, $\alpha^1$ and $\alpha^2$ 
had a common tail $c\in\kQ_0\setminus\{p^0,\dots,p^l\}$ of valence four,
then
the two arrows having $c$ as common head could not occur in paths avoiding $b$
and leading from $p^0$ to $\{p^1,\dots,p^l\}$.   
In particular, $\kQP=\kQ\setminus\{b\}$ would not be stable.
The reversed implication may be shown along the lines of the proof of
Lemma \ref{lem-polytrees}. 
\vspace{1ex}\\
It remains to deal with the polytope $\kDstQ$ itself.
If $\dim\kDstQ\geq 4$, then we have just shown that every
facet is Minkowski prime; this implies the same property for $\kDstQ$, too.
The three-dimensional case can be easily solved by a complete classification.
\end{proof}

\begin{corollary}\klabel{cor-T1}
For $X_\kDstQ$, the non-trivial $T_X^1(-R)$ are
one-dimensional, and they arise from degrees $R=R(b)$ such that
$R(b)\leq 1$ on $\kDstQ$ and $R(b)=1$ exactly on 
$\kPF(\kQ\setminus b)$ with $b$ being a simple knot.
\end{corollary}

\neu{deform-T1+}
The precise description of $T_X^1$ given in Corollary \ref{cor-T1}
may be supplemented by the following, more structural claim.
\par

\begin{theorem}\klabel{th-omega}
Let $X_\kDstQ$ be the toric Gorenstein singularity assigned to the quiver
polytope of a flag quiver $\kQ$. Then, the simple knots
\vspace{-0.5ex}
$\,b\in\kQ_0\setminus\{p^0,\dots,p^l\}$ 
parametrize the 3-codimensional $A_1$-strata 
$
\xymatrix@!0@C=4.5em@R=0ex{
Z(b)\,\ar@{^(->}[r]^-{i(b)}& X_\kDstQ
}
$,
and the module of the infinitesimal deformations of $X_\kDstQ$ equals
\[
T^1_X
={\ks\bigoplus}_b \,i(b)_\ast \;\omega_{Z(b)} 
\otimes \omega_{X}^{-1}
\]
with $\omega_{\ldots}$ denoting the canonical sheaves on $Z(b)$ and $X$.
\end{theorem}

\begin{proof}
We know from Corollary \ref{cor-two} that the 
3-codimensional singularities in $X_\kDstQ$ correspond to the 2-dimensional, 
non-triangular faces of $\kDstQ$ which are squares. On the other hand,
Proposition \ref{prop-simpknot} establishes their relation to the simple 
knots $b$. For any of it we may define 
\[
T(b)\;:=\;\oplus_{R(b)} \,T^1\big(-R(b)\big)
\;=\;\kT^1\big(\kPF(\kQ\setminus b)
\big)^{\# \{R(b)'\mbox{\kf s}\}}
\]
meaning to sum over all $R(-b)$ in the sense of Corollary \ref{cor-T1}.
Thus, the whole $T^1_X$ splits, as a $\C$-vector space, into
$T^1_X = \oplus_b\, T(b)$.
\vspace{1.0ex}\\
The module structure of $T^1_X$ has been explained in 
Theorem \ref{th-TD}: 
The dual cone of $\sigma=\cone(\kDstQ)\subseteq N_\R$
is $\sigma^\vee=\cone(\kDQ)\subseteq M_\R$ 
with $M:=\gHom(N,\Z)=\kF\oplus\Z$, cf.\ \zitat{quiver}{weights}.
Whenever $s\in \sigma^\veee\cap M$ vanishes
on $\kPF(\kQ\setminus b)$, i.e.\ if 
$s\in \cone\,\kD(\kQ\setminus b,\, \wtc_\kQ)\subseteq\cone\,\kDQ$,
then the multiplication $x^s:T^1\big(-R(b)\big)\to T^1\big(-R(b)+s\big)$ 
is the identity map
when both sides are identified with the one-dimensional
$\kT^1\big(\kPF(\kQ\setminus b)\big)$.
If $s$ does not vanish on $\kPF(\kQ\setminus b)$, 
then the multiplication is zero.\\
Hence, the splitting of $T^1$ respects the module structure. Moreover,
on the summands $T(b)$, this structure factors through the surjection
$\C[\sigma^\veee\cap M]\surj 
\C[\sigma^\veee\cap \kPF(\kQ\setminus b)^\bot\cap M]=
\C[\cone\kD(\kQ\setminus b,\, \wtc_\kQ)\cap M]$ with
\[
T(b)=\Big(\oplus_{R\in\,\mbox{\kf int}\,
\mbox{\kf cone}\,\kD(\kQ\setminus b,\, \wtc_\kQ)} 
\,\C\cdot x^{R-\ku{1}}\Big) \otimes_\C \kT^1\big(\kPF(\kQ\setminus b)\big)\,.
\]
On the other hand, the semigroup algebra 
$\C[\cone\kD(\kQ\setminus b,\, \wtc_\kQ)\cap M]$ equals the coordinate ring
of the stratum $Z(b)$, and it is a general fact for affine toric varieties
$\toric{\tau} = \Spec\, \C[\tau^\vee\cap M]$ that 
$\oplus_{R\in\,\mbox{\kf int}\,\tau^\vee}\C\cdot x^R$ equals the canonical
module $\omega_{\toric{\tau}}$.
\end{proof}

\neu{deform-height0}
Eventually, we would rather like to study the deformations of the projective
toric variety $\PP_\kDQ$ instead of that of its cone $X_\kDstQ$. 
This just means to focus on those multidegrees $R$ with height
or $\Z$-degree $0$,
i.e.\ on $\,R\in\kF\times\{0\}\subseteq \kF\times\Z=M$.
To use Corollary \ref{cor-T1} for describing the entire
homogeneous piece $T^1_X(0)$,
it is helpful to realize $M$ as a subspace of $\Z^{\kQ_1}$.
This is done by the isomorphism
\[
\xymatrix@!0@C=6.5em@R=3.5ex{
\kF\oplus \Z \ar[r]^-{\sim} & \pi^{-1}(\Z\cdot\wtc)\,,&
(r,g) \ar@{|->}[r] & r + g\,\ku{1}\\
&& (\kD-\ku{1},\,1) \ar@{|->}[r] & \kD\,.
}
\]
Under this map,
the value of $R=(r,g)\in M$ on the
vertex $(a^\alpha,1)$ of $(\kD-\ku{1},\,1)$ 
equals the $[\alpha]$-coordinate of $R=r + g\,\ku{1}\in\Z^{\kQ_1}$.
In particular, multidegrees $R$ of height $0$
are exactly those coming from
$\pi^{-1}(0\cdot\wtc)=\kF\subseteq\Z^{\kQ_1}$.

\begin{definition}\klabel{def-multipath}
Let $b$ be a simple knot and denote by $a^1, a^2$ the tails of the two arrows
$\alpha^1,\alpha^2$ 
with head $b$, respectively; 
$c^1,c^2$ are defined in a similar manner on the outgoing arrows
$\beta^1,\beta^2$.
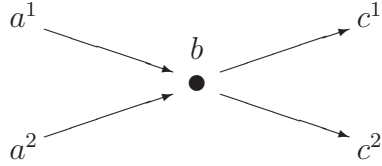
\begin{center}
\unitlength=0.6pt
\begin{picture}(203.00,83.00)(51.00,738.00)
\put(41.00,812.00){\vector(3,-1){81.00}}
\put(152.00,785.00){\vector(3,1){81.00}}
\put(41.00,744.00){\vector(3,1){81.00}}
\put(152.00,771.00){\vector(3,-1){81.00}}
\put(137.00,778.00){\circle*{10}}
\put(28.00,821.00){\makebox(0,0)[cc]{$a^1$}}
\put(28.00,738.00){\makebox(0,0)[cc]{$a^2$}}
\put(246.00,821.00){\makebox(0,0)[cc]{$c^1$}}
\put(246.00,738.00){\makebox(0,0)[cc]{$c^2$}}
\put(137,800){\makebox(0,0)[cc]{$b$}}
\end{picture}
\end{center}
An element $R\in\Z^{\kQ_1}_{\geq 0}$ is called {\em multipath}
from $\{a^1,a^2\}$ to $\{c^1,c^2\}$
if $\,\pi(R)=[a^1]+[a^2]-[c^1]-[c^2]$.
The standard example is $R_b:=[\alpha^1]+[\alpha^2]+[\beta^1]+[\beta^2]$
through $b$.
\vspace{1.5ex}
\end{definition}

\begin{theorem}\klabel{th-multipath}
Running through the simple knots $\,b\in\kQ_0\setminus\{p^0,\dots,p^l\}$,
the part
$T^1_X(0)$ splits into
$T^1_X(0)=\oplus_b T_0(b)$, and the dimension of each $T_0(b)$
equals the number of multipaths
leading from $\{a^1,a^2\}$ to $\{c^1,c^2\}$, but avoiding $b$.
\vspace{-1.5ex}
\end{theorem}

\begin{proof}
Corollary \ref{cor-T1} characterizes the $T^1$-carrying multidegrees 
$R(b)\in\Z^{\kQ_1}$ assigned to the knot $b$ by the properties
\kitem
\item
$R(b)_\alpha=1$ for $\alpha$ being one of the four arrows touching $b$
and
\vspace{0.5ex}
\item
$R(b)_\alpha\leq 0$ for the remaining arrows $\alpha\in(\kQ\setminus b)_1$.
\kenditem
On the other hand, the condition of having height $0$
means $R(b)\in\kF$, i.e.\ that $R(b)$ encodes an
cycle inside $\kQ$. 
Hence, the negative multidegrees $-R(b)$ 
represent cycles 
using each of the four $b$-arrows exactly once and in the wrong direction,
but respecting the orientation of the remaining arrows in $\kQ$.
With other words, $R_b-R(b)$ represents a multipath
from $\{a^1,a^2\}$ to $\{c^1,c^2\}$ avoiding $b$.
\end{proof}

{\bf Example \ref{ex-fence}} (continued){\bf .}
While, in the quiver $\kQ^1$ of \zitat{flag}{poly}, the vertices $a$ and $b$
give rise to unique multidegrees $R(a)$ and $R(b)$ of height $0$, there
are five multipaths 
corresponding to $c$. Leaving out $R_c$, the remaining four paths do not touch
$c$. Hence, they are responsible for four dimensions 
inside the six-dimensional $T^1_X(0)$.
\par

\neu{deform-noheight0}
Whenever $\kPF<\kDstQ$ is a face, then there exist always degrees 
$R\in M=\pi^{-1}(\Z\cdot\wtc)$ such that $R\leq 1$ on $\kDstQ$ and
$R=1$ exactly on $\kPF$ -- just take $R$ as the difference of $\ku{1}$ and
an interior lattice point of the $\sigma^\veee$-face dual to $\kPF$.
In particular, as we have already seen in Theorem \ref{th-omega},
every simple vertex $b$ provides a contribution to $T^1_X$.\\
However, it might happen that simple vertices $b$ do not provide multidegrees
$R(b)$ of height $0$. In the following  
{\em non-plane} flag quiver $Q^3$,
every of the five inner vertices is simple in the sense of       
Definition \ref{def-simpknot}. While the knots $c^1,\dots,c^4$ provide
multidegrees $R(c^i)$ of height $0$, the knot $b$ does not. The reason
is that there is exactly one multipath leading from $\{c^3,c^4\}$ to
$\{c^1,c^2\}$, but this multipath touches $b$.
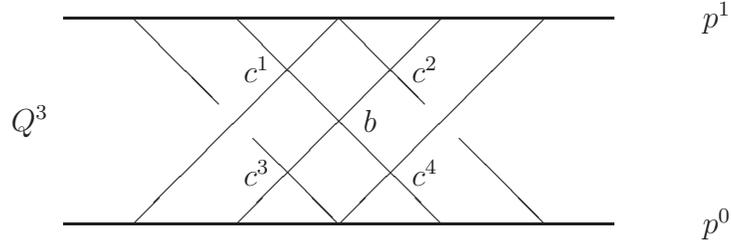
\begin{center}
\hspace*{\fill}
\unitlength=1.3pt
\begin{picture}(180.00,60.00)(40.00,740.00)
\thicklines
\put(40.00,800.00){\line(1,0){160.00}}
\put(40.00,740.00){\line(1,0){160.00}}
\thinlines
\put(120.00,800.00){\line(-1,-1){60}}
\put(150.00,800.00){\line(-1,-1){60}}
\put(180.00,800.00){\line(-1,-1){60}}
\put(60.00,800.00){\line(1,-1){25}}
\put(95.00,765.00){\line(1,-1){25}}
\put(90.00,800.00){\line(1,-1){60}}
\put(120.00,800.00){\line(1,-1){25}}
\put(155.00,765.00){\line(1,-1){25}}
\put(30.00,770.00){\makebox(0,0)[cc]{$\kQ^3$}}
\put(96.00,785.00){\makebox(0,0)[cc]{$c^1$}}
\put(145.00,785.00){\makebox(0,0)[cc]{$c^2$}}
\put(96.00,755.00){\makebox(0,0)[cc]{$c^3$}}
\put(145.00,755.00){\makebox(0,0)[cc]{$c^4$}}
\put(129.00,770.00){\makebox(0,0)[cc]{$b$}}
\put(230.00,740.00){\makebox(0,0)[cc]{$p^0$}}
\put(230.00,800.00){\makebox(0,0)[cc]{$p^1$}}
\end{picture}
\hspace*{\fill}
\end{center}
Here is another example. The {\em plane} flag quiver $Q^4$ is built from
$Q^1$ of Example \ref{ex-fence} in \zitat{flag}{poly} by 
adding one single rope. However, this procedure removes the
$b$-contribution from $T^1_X(0)$,
i.e.\ $T_0(b)=0$.
\vspace{-4ex}
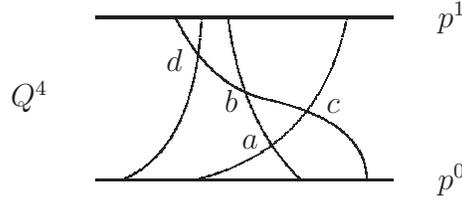
\begin{center}
\unitlength=0.5pt
\begin{picture}(220.00,140.00)(-20.00,680.00)
\thicklines
\put(00.00,800.00){\line(1,0){225.00}}
\put(00.00,677.00){\line(1,0){225.00}}
\thinlines
\bezier106(60.00,800.00)(86.00,746.00)(133.00,737.00)
\bezier106(133.00,737.00)(204.00,720.00)(205.00,677.00)
\bezier106(100.00,800.00)(110.00,720.00)(155.00,677.00)
\bezier106(75.00,677.00)(170.00,700.00)(190.00,800.00)
\bezier106(20.00,677.00)(75.00,700.00)(80.00,800.00)
\put(60.00,768.00){\makebox(0,0)[cc]{$d$}}
\put(-50.00,740.00){\makebox(0,0)[cc]{$\kQ^4$}}
\put(117.00,706.00){\makebox(0,0)[cc]{$a$}}
\put(180.00,733.00){\makebox(0,0)[cc]{$c$}}
\put(103.00,737.00){\makebox(0,0)[cc]{$b$}}
\put(270.00,677.00){\makebox(0,0)[cc]{$p^0$}}
\put(270.00,800.00){\makebox(0,0)[cc]{$p^1$}}
\end{picture}
\end{center}
In general,
the absence of $T_X^1(0)$-pieces for a simple knot of $\kQ$
means that the corresponding $A_1$-singularity is, even locally,
not smoothable with deformations of degree $0$. Hence,
to obtain smoothability of $\PP_{\kDQ}$ in codimension three,
a neccessary condition is that $T_0(b)\neq 0$ for every simple knot
$b$ of $\kQ$. Using Theorem \ref{th-multipath},
this translates into the existence of detouring multipaths
around every simple knot.

\neu{deform-T2}
We will show that the just mentioned neccessary condition
for smoothability in codimension three is sufficient, too.
For simple knots $b$,
the one-parameter deformations of $X_\kDstQ$
provided by the one-dimensional vector
spaces $T^1_X\big(-R(b)\big)$ are smoothings of the 
three-codimensional $\mbox{\rm A}_1$-singularities along $Z(b)$.
The latter are the orbits of the
three-dimensional cones over the faces 
$\kPF(\kQ\setminus b)\leq\kDstQ$.
The question whether these one-parameter
families fit together in a huge deformation doing all the smoothings
simultaneously leads to the investigation of $T^2_X$.

\begin{theorem}\klabel{th-lackobst}
Whenever $R\in M$ is a positive linear combination of degrees
$R^i\in M$ carrying infinitesimal deformations, i.e.\
$T^1_X(-R^i)\neq 0$, then $T^2_X(-R)=0$.
\end{theorem}

\begin{proof}
Because of Corollary \ref{cor-two} and
Theorem \ref{th-TD}, we may assume that $R\leq 1$ on $\kDstQ$. 
To apply Theorem \ref{th-D2van}, we have first
to check the three-dimensional faces 
of $\kDstQ\cap [R=1]\subseteq\kDstQ$ for
non-pyramids, i.e.\  
to exclude octahedra and prisms corresponding to the
triangular quivers depicted in \zitat{tight}{smallfaces}. While the latter
would provide a three-dimensional face with non-trivial $\kT^1$, which is
excluded by Proposition \ref{prop-simpknot}, we need a finer argument for the
octahedra:
\vspace{1ex}\\
The assumption of our theorem says that 
$R^i_\alpha=\langle a^\alpha, R^i\rangle \geq 1$
(in fact ``$\,=\!1$'') holds exactly
for the four arrows $\alpha$ containing the simple knot $b(R^i)$,
cf.\ Corollary \ref{cor-T1} or the proof of Theorem \ref{th-multipath}.
In particular, since $R$ is a positive linear combination of
those $R^i$, the relation $\langle a^\alpha, R\rangle \geq 1$ is 
impossible, unless $\tail(\alpha)$ or $\head(\alpha)$ is a simple vertex.
On the other hand, if $R=1$ was true on an octahedron $\kPF(\kQP)$, i.e.\
on the arrows of $\kG_\kQ(\kQP)\cong\kGdbl(3)$,
then two of the three vertices of $\kG_\kQ(\kQP)$ would equal 
original vertices from 
$\kQ_0\setminus\{p^0,\ldots,p^l\}$.
However,
as for $\kQ^2$ in Example \ref{ex-fence} of \zitat{deform}{T1}, 
the double arrow between these vertices implies that they
cannot be simple -- providing a contradiction.
\vspace{1ex}\\
Let us now check the remaining assumptions
of Theorem \ref{th-D2van}. The
two-dimensional, non-triangular faces of $\kDstQ\cap [R=1]$ are squares
provided by simple knots $b\in\kQ_0$; the four vertices of these squares
correspond to the arrows containing $b$. In particular, the property of an
arrow to contain exactly two knots translates into the property of an
vertex $a^\alpha$ of $\kDstQ$ to belong to at most two of these squares.
This means that we are done in case of
$\dim \big(\kDstQ\cap [R=1]\big)\geq 5$.\\
For the remaining case $\dim \big(\kDstQ\cap [R=1]\big)=4$, our argument
requires a slight refinement. To obtain vanishing of
$\kT^2$, we do not use 
Theorem \ref{th-D2van} itself, 
but the stronger, original Theorem (4.7)
of \cite{hodge}: Since the quiver $\kQ$ lacks oriented cycles,
it provides a (non-linear) ordering of the set $\kQ_0$. 
Hence,
whenever there is a connected set of squares in
our face $[R=1]$, then there is at least one vertex of one of these squares
that contains only this single square. Now, beginning with this particular
vertex, we may ``clean'' these squares in the sense of  \cite{hodge}, (4.7)
successively.
\end{proof}

\begin{remark}
It is not true in general that $T^2_X=0$, cf.\
Example \ref{ex-oct}.
However, the previous theorem says that at least the obstructions inside
$T^2_X$ are void.
\vspace{1ex}
\end{remark}

\begin{corollary}\klabel{cor-smoothable}
Gorenstein singularities provided by flag quivers are smoothable in
codimension three. Moreover, if every simple knot $b$ can be by-passed
with a multipath
connecting its neighbors,
then this can be done by a deformation of degree $0$.
\vspace{-1ex}
\end{corollary}

\begin{proof}
With the notation of Corollary \ref{cor-T1}, 
we choose one element $R(b)\in M$ for each
simple knot $b$. By the lack of obstructions, the corresponding one-parameter
families fit into a common deformation over a smooth parameter space $S$.\\
Now, looking at the general points of the singular three-codimensional
strata, $S$ is obtained from their one-dimensional versal deformations via
base change. In particular, for each of these strata, there is a hypersurface
in $S$ containing the parameters {\em not} smoothing this stratum.
Hence, taking a curve inside $S$ that avoids all these hypersurfaces outside
$0\in S$, yields the desired smoothing.
\vspace{1ex}
\end{proof}

\begin{example}\klabel{ex-smoothable}
The $5$-dimensional projective varieties $\PP_\kDQ$ corresponing to the
quivers $\kQ^1$ iand $\kQ^2$ of \zitat{flag}{poly} 
are smoothable in codimension three.
On the other hand, for the quivers 
$\kQ^3$ and $\kQ^4$ of \zitat{deform}{noheight0}, we
know this only for 
$X_{\kDstQ}=\Cone(\PP_\kDQ)$
instead for $\PP_\kDQ$ itself.
\end{example}

{\small
\parbox{7cm}{
Klaus Altmann\\
Fachbereich Mathematik und Informatik\\
Freie Universit\"at Berlin\\
Arnimalle 3\\
14195 Berlin, Germany\\
e-mail: altmann@math.fu-berlin.de}
\hfill
\parbox[r]{8cm}{
\hspace*{\fill}Duco van Straten\\
\hspace*{\fill}Fachbereich Mathematik (17)\\
\hspace*{\fill}Johannes Gutenberg-Universit\"at\\
\hspace*{\fill}D-55099 Mainz, Germany%
\vspace{1ex}\\
\hspace*{\fill}e-mail: straten@mathematik.uni-mainz.de}}

\end{document}